\documentclass[10pt]{article}
\usepackage{amsmath,amsthm,amsfonts,amssymb,bbm,
latexsym,cite}
\usepackage{psfrag,graphicx}

\newtheorem{lemme}{Lemma}[section]
\newcommand{\bl}{\begin{lemme}}
\newcommand{\el}{\end{lemme}}

\newcommand{\brem}{\begin{remark}}
\newcommand{\erem}{\end{remark}}
\newcommand{\bconj}{\begin{conjecture}}
\newcommand{\econj}{\end{conjecture}}
\newcommand{\bdefi}{\begin{definition}}
\newcommand{\edefi}{\end{definition}}
\newcommand{\bt}{\begin{theo}}
\newcommand{\bfa}{\begin{fact}}
\newcommand{\efa}{\end{fact}}

\newcommand{\et}{\end{theo}}
\newcommand{\bp}{\begin{prop}}
\newcommand{\ep}{\end{prop}}
\newcommand{\be}{\begin{equation}}

\newcommand{\ee}{\end{equation}}

\usepackage{graphicx}
 \usepackage{mathptmx}      
\usepackage[OT1]{fontenc}
\usepackage{amsmath}
\usepackage[colorlinks,citecolor=blue,urlcolor=blue]{hyperref}
\usepackage{hypernat}
\usepackage{latexsym,graphicx,amssymb}

\newtheorem{theo}{Theorem}[section]

\newcommand{\E}{{\mathbb{E}}}
\newcommand{\supp}{{\sf Supp}}
\newcommand{\im}{{\sf Im}}
\newcommand{\re}{{\sf Re}}
\newcommand{\tr}{{\sf Tr}}
\newcommand{\diag}{{\sf Diag}}
\newcommand{\ith}{{i^{\rm th}}}

\newcommand{\kth}{{k^{\rm th}}}

\newcommand{\bone}{{\bf 1}}
\newcommand{\C}{{\mathbb{C}}}
\newcommand{\R}{{\mathbb{R}}}
\newtheorem{rem}{Remark}
\begin{document}

\title{Eigenvectors of some large sample covariance matrix ensembles.}



\author{Olivier Ledoit  and
        Sandrine P\'ech\'e \\ \\
              Institute for Empirical Research in Economics, University of Zurich,\\
 Bl\"{u}mlisalpstrasse 10, 8006 Z\"{u}rich, Switzerland\\ 
              oledoit@iew.uzh.ch \\ \\          
              Institut Fourier, Universit\'e Grenoble 1,\\
 100 rue des Maths, BP 74, 38402 Saint-Martin-d'H\`{e}res, France\\
sandrine.peche@ujf-grenoble.fr
}

\date{}

\maketitle
 
\begin{abstract}
We consider sample covariance matrices $S_N=\frac{1}{p}\Sigma_N^{1/2}X_NX_N^* \Sigma_N^{1/2}$ where $X_N$ is a $N \times p$ real or complex matrix with i.i.d. entries with finite $12^{\rm th}$ moment and $\Sigma_N$ is a $N \times N$ positive definite matrix. In addition we assume that the spectral measure of $\Sigma_N$ almost surely converges to some limiting probability distribution as $N \to  \infty$ and $p/N \to \gamma >0.$
We quantify the relationship between sample and population eigenvectors
by studying the asymptotics of functionals of the type $\frac{1}{N} \text{Tr} \left ( g(\Sigma_N) (S_N-zI)^{-1})\right ),$
where $I$ is the identity matrix, $g$ is a bounded function and $z$ is a complex number. 
This is then used to compute the asymptotically optimal bias correction for sample eigenvalues, paving the way for a new generation of improved estimators of the covariance matrix and its inverse.
 
\end{abstract}

\section{Introduction and Overview of the Main Results}
\label{sec:introduction}
\subsection{Model and results \label{subsec:model}}
Consider $p$ independent samples $C_1,\ldots,$ $C_p$, all of which are $N \times 1$ real or complex vectors. 
In this paper, we are interested in the large-$N$-limiting spectral properties of the sample covariance matrix 
$$S_N= \frac{1}{p} CC^*, \quad C=[C_1, C_2, \ldots , C_p],$$
when we assume that the sample size $p=p(N)$ satisfies $p/N \to \gamma $ as $N \to \infty$ for some $\gamma>0.$ 
This framework is known as large-dimensional asymptotics. 
Throughout the paper, $\bone$ denotes the indicator function of a set, and we make the following assumptions: $C=\Sigma_N^{1/2}X_N$ where 
\begin{itemize}
\item{$(H_1)$} $X_N$ is a $N\times p$ matrix of real or complex iid random variables with zero mean, unit variance, and $12^{\rm th}$ absolute central moment bounded by a constant $B$ independent of $N$ and $p$;
\item{$(H_2)$} the population covariance matrix $\Sigma_N$ is a $N$-dimensional random Hermitian positive definite matrix independent of $X_N$;
\item{$(H_3)$} $p/N\to \gamma>0$ as $N\to\infty$;
\item{$(H_4)$} $(\tau_1,\ldots,\tau_N)$ is a system of eigenvalues of $\Sigma_N$, and the empirical spectral distribution  (e.s.d.) of the population covariance given by $H_N(\tau)=\frac{1}{N}\sum_{j=1}^N\bone_{[\tau_j,+\infty)}(\tau)$ converges a.s. to a nonrandom limit $H(\tau)$ at every point of continuity of $H$. $H$ defines a probability distribution function, whose support $\supp(H)$ is included in the compact interval $[h_1,h_2]$ with $0< h_1\leq h_2<\infty$.
\end{itemize}
The aim of this paper is to investigate the asymptotic properties of the eigenvectors of such sample covariance matrices. In particular, we will quantify how the eigenvectors of the sample covariance matrix deviate from those of the population covariance matrix under large-dimensional asymptotics. This will enable us to characterize how the sample covariance matrix deviates {\em as a whole} (i.e.~through its eigenvalues {\em and} its eigenvectors) from the population covariance matrix. Specifically, we will introduce bias-correction formulae for the eigenvalues of the sample covariance matrix that can lead, in future research, to improved estimators of the covariance matrix and its inverse. This will be developped in the discussion (Sections \ref{sec:eigenvectors} and \ref{sec:shrinkage}) following our main result Theorem \protect\ref{theo:theta} stated below.\\ 
Before exposing our results, we briefly review some known results about the spectral properties of sample covariance matrices under large-dimensional asymptotics.\\
In the whole paper we denote by $((\lambda_1^N,\ldots,\lambda_N^N);(u_1^N,\ldots,u_N^N))$ a system of eigenvalues and orthonormal eigenvectors of the sample covariance matrix $S_N=\frac{1}{p}\Sigma_N^{\frac{1}{2}}X_N^{}X_N^*\Sigma_N^{\frac{1}{2}}$.\label{item:assumption}
Without loss of generality, we assume that the eigenvalues are sorted in decreasing order: $\lambda_1^N\geq\lambda_2^N\geq\cdots\geq\lambda_N^N$. We also denote by $(v_1^N,\ldots,v_N^N)$ a system of orthonormal eigenvectors of $\Sigma_N$. Superscripts will be omitted when no confusion is possible.

\paragraph{} First the asymptotic behavior of the eigenvalues is now quite well understood. The ``global behavior'' of the spectrum of $S_N$ for instance is characterized through the e.s.d., defined as: $ F_N(\lambda)=N^{-1}\sum_{i=1}^N\bone_{[\lambda_i,+\infty)}(\lambda)$, $\forall \lambda \in \mathbb{R}$. The e.s.d.~is usually described through its Stieltjes transform. We recall that the Stieltjes transform of a nondecreasing function $G$ is defined by $m_G(z)=\int_{-\infty}^{+\infty}(\lambda-z)^{-1}dG(\lambda)$ for all $z$ in $\C^+$, where $\mathbb{C}^+=\{z \in \mathbb{C},\, \im(z)>0\}.$
The use of the Stieltjes transform is motivated by the following inversion formula: given any nondecreasing function $G$, one has that  $G(b)-G(a)=\lim_{\eta\to0^+}\pi^{-1}\int_a^b\im\left[m_G(\xi+i\eta)\right]d\xi$, which holds if $G$ is continuous at $a$ and $b$.

The first fundamental result concerning the asymptotic global behavior of the spectrum has been obtained by Mar\v{c}enko and Pastur in \cite{MP}. Their result has been later precised e.g. in \cite{SilversteinBai,GS77,YK83,W78,Y86}. 
In the next Theorem, we recall their result (which was actually proved in a more general setting than that exposed here) and quote the most recent version as given in \cite{Silverstein}. 

Let $m_{F_N}(z)=\frac{1}{N}\sum_{i=1}^N\frac{1}{\lambda_i-z}
=\frac{1}{N}\tr\left[(S_N-zI)^{-1}\right],$
where $I$ denotes the $N \times N$ identity matrix.

\begin{theo}[\cite{MP}]
Under Assumptions $(H_1)$ to $(H_4)$, one has that for all $z \in \C^+,$ \\$\lim_{N \to \infty }m_{F_N}(z)=m_F(z)$ a.s. where
\be
\label{eq:marcenko-pastur}
\forall z\in\mathbb{C}^+,\quad m_F(z)=\int_{-\infty}^{+\infty}\left\{\tau\left[1-\gamma^{-1}-\gamma^{-1}z\,m_F(z)\right]-z\right\}^{-1}dH(\tau).
\ee
Furthermore, the e.s.d.~of the sample covariance matrix given by $F_N(\lambda)=N^{-1}\sum_{i=1}^N
\bone_{[\lambda_i,+\infty)}(\lambda)$ converges a.s.~to the nonrandom limit $F(\lambda)$ at all points of continuity of $F$.
\end{theo} 
In addition, \cite{SC95} show that the following limit exists :
\be \forall\lambda\in\R-\{0\},\, \lim_{z\in\mathbb{C}^+\to\lambda}m_F(z)\equiv \breve{m}_F(\lambda)\label{defbrevem}.\ee They also prove that $F$ has a continuous derivative which is given by $F'=\pi^{-1}\im[\breve{m}_F]$ on $(0,+\infty)$. 
More precisely, when $\gamma>1,$ $\lim_{z\in\mathbb{C}^+\to\lambda}m_F(z)\equiv \breve{m}_F(\lambda)$ exists for \emph{all} $\lambda\in\mathbb{R}$, $F$ has a continuous derivative $F'$ on \emph{all} of $\mathbb{R}$, and $F(\lambda)$ is identically equal to zero in a neighborhood of $\lambda=0$.
When $\gamma<1$, the proportion of sample eigenvalues equal to zero is asymptotically $1-\gamma$. In this case, 
it is convenient to introduce the e.s.d.~$\underline{F}=\left(1-\gamma^{-1}\right)\bone_{[0,+\infty)}+\gamma^{-1}F$,
which is the limit of e.s.d.~of the $p$-dimensional matrix $p^{-1}X_N^*\Sigma_NX_N$. 
Then $\lim_{z\in\mathbb{C}^+\to\lambda}m_{\underline{F}}(z)\equiv \breve{m}_{\underline{F}}(\lambda)$ exists for \emph{all} $\lambda\in\mathbb{R}$, $\underline{F}$ has a continuous derivative $\underline{F}'$ on \emph{all} of $\mathbb{R}$, and $\underline{F}(\lambda)$ is identically equal to zero in a neighborhood of $\lambda=0$. 
When $\gamma$ is exactly equal to one, further complications arise because the density of sample eigenvalues can be unbounded in a neighborhood of zero; for this reason we will sometimes have to rule out the possibility that $\gamma=1$.\\
Further studies have complemented the a.s. convergence established by the Mar\v{c}enko-Pastur theorem (see e.g. \cite{B93,BS98,BS99,BS04,BY93,J01,Peche} and \cite{BaiLarge} for more references).  
The Mar\v{c}enko-Pastur equation has also generated a considerable amount of interest in statistics \cite{EK08,LW02}, finance \cite{L95,Honey}, signal processing \cite{SC92}, and other disciplines. We refer the interested reader to the recent book by Bai and Silverstein \cite{BS06} for a throrough survey of this fast-growing field of research.

\paragraph{}
As we can gather from this brief review of the literature, the Mar\v{c}en\-ko-Pastur equation reveals much of the behavior of the \emph{eigenvalues} of sample covariance matrices under large-dimensional asymptotics. 
It is also of utmost interest to describe the asymptotic behavior of the \emph{eigenvectors}. Such an issue is fundamental to statistics (for instance both eigenvalues and eigenvectors are of interest in Principal Components Analysis), communication theory (see e.g. \cite{PanZhou08}), wireless communication, finance. The reader is referred to \cite{BMP07}, Section 1 for more detail and to \cite{BickelLevina} for a statistical approach to the problem and a detailed exposition of statistical applications.
\\
Actually much less is known about eigenvectors of sample covariance matrices. In the special case where $\Sigma= I$ and the $X_{ij}$ are i.i.d. standard (real or complex) Gaussian random variables, it is well-known that the matrix of sample eigenvectors is Haar distributed (on the orthogonal or unitary group). To our knowledge, these are the only ensembles for which the distribution of the eigenvectors is explicitly known.
It has been conjectured that for a wide class of non Gaussian ensembles, the matrix of sample eigenvectors should be ``asymptotically Haar distributed'', provided $\Sigma=I$. Note that the notion  ``asymptotically Haar distributed'' needs to be defined. 
This question has been investigated by \cite{S84}, \cite{S89}, \cite{S90} followed by \cite{BMP07} and \cite{PanZhou08}. Therein a random matrix $U$ is said to be asymptotically Haar distributed if $Ux$ is asymptotically uniformly distributed on the unit sphere for any non random unit vector $x$. \cite{S90} and \cite{BMP07} are then able to prove the conjecture under various sets of assumptions on the $X_{ij}$'s.\\
In the case where $\Sigma \not=I$, much less is known (see \cite{BMP07} and \cite{PanZhou08}). One expects that the distribution of the eigenvectors is far from being rotation-invariant. This is precisely the aspect in which this paper is concerned. 
\paragraph{}
In this paper, we present another approach to study eigenvectors of sample covariance matrices. Roughly speaking, we study ``functionals'' of the type 
\begin{eqnarray}
\forall z\in\mathbb{C}^+,\qquad
\Theta^g_N(z)&=&\frac{1}{N}\sum_{i=1}^N\frac{1}{\lambda_i-z}\sum_{j=1}^N\left|u_i^*v_j\right|^2\times g(\tau_j)\label{eq:general}\\
&=&\frac{1}{N}\tr\left[(S_N-zI)^{-1}g(\Sigma_N)\right],\nonumber
\end{eqnarray}
where $g$ is any real-valued univariate function satisfying suitable regularity conditions. By convention, $g(\Sigma_N)$ is the matrix with the same eigenvectors as $\Sigma_N$ and with eigenvalues  $g(\tau_1),\ldots,g(\tau_N)$.
These functionals are generalizations of the Stieltjes transform used in the Mar\v{c}enko-Pastur equation.
Indeed, one can rewrite the Stieltjes transform of the e.s.d.~as:
\be
\label{eq:mFN}
\forall z\in\mathbb{C}^+,\qquad
m_{F_N}(z)=\frac{1}{N}\sum_{i=1}^N\frac{1}{\lambda_i-z}\sum_{j=1}^N\left|u_i^*v_j\right|^2\times1.
\ee
The constant $1$ that appears at the end of Equation (\ref{eq:mFN}) can be interpreted as a weighting scheme placed on the population eigenvectors: specifically, it represents a flat weighting scheme.
The generalization we here introduce puts the spotlight on how the sample covariance matrix relates to the population covariance matrix, or even \emph{any function of} the population covariance matrix.

Our main result is given in the following Theorem. 
\bt
\label{theo:theta}
Assume that conditions $(H_1)-(H_4)$ are satisfied.
Let $g$ be a (real-valued) bounded function defined on $[h_1,h_2]$ with finitely many points of discontinuity. Then there exists a nonrandom function $\Theta^g$ defined over $\mathbb{C}^+$ such that $\Theta_N^g(z)=N^{-1}\tr\left[(S_N-zI)^{-1}g(\Sigma_N)\right]$ converges a.s.~to $\Theta^g(z)$ for all $z\in\mathbb{C}^+$. Furthermore, $\Theta^g$ is given by:
\be
\label{eq:theta}
\forall z\in\mathbb{C}^+,\;\Theta^g(z)
=\int_{-\infty}^{+\infty}\left\{\tau\left[1-\gamma^{-1}-\gamma^{-1}zm_F(z)\right]-z\right\}^{-1}g(\tau)dH(\tau).
\ee
\et

One can first observe that as we move from a flat weighting scheme of $g\equiv 1$ to any arbitrary weighting scheme $g(\tau_j)$, the integration kernel $\left\{\tau\left[1-\frac{1}{\gamma}-\frac{z\,m_F(z)}{\gamma}\right]-z\right\}^{-1}$ remains unchanged. Therefore, our Equation (\ref{eq:theta}) generalizes Mar\v{c}enko and Pastur's foundational result. 
Actually the proof of Theorem \ref{theo:theta} follows from some of the arguments used in \cite{Silverstein} to derive the Marchenko-Pastur equation. This proof is postponed until Section \ref{sec:prooftheta}.

The generalization of the Mar\v{c}enko-Pastur equation we propose allows to consider a few unsolved problems regarding the overall relationship between sample and population covariance matrices. 
Let us consider two of these problems, which are investigated in more detail in the two next subsections.\\
The first of these questions is: how do the eigenvectors of the sample covariance matrix deviate from those of the population covariance matrix? By injecting functions $g$ of the form $\bone_{(-\infty,\tau)}$ into Equation (\ref{eq:theta}), we quantify the asymptotic relationship between sample and population eigenvectors. This is developed in more detail in Section \ref{sec:eigenvectors}.\\
Another question is: how does the sample covariance matrix deviate from the population covariance matrix as a whole, and how can we modify it to bring it closer to the population covariance matrix? This is an important question in Statistics, where a covariance matrix estimator that improves upon the sample covariance matrix is sought. By injecting the function $g(\tau)=\tau$ into Equation (\ref{eq:theta}), we find the optimal asymptotic bias correction for the eigenvalues of the sample covariance matrix in Section \ref{sec:shrinkage}. We also perform the same calculation for the \emph{inverse} covariance matrix (an object of great interest in Econometrics and Finance), this time by taking $g(\tau)=1/\tau$.\\
This list is not intended to be exhaustive. Other applications may hopefully be extracted from our generalized Mar\v{c}enko-Pastur equation.

\subsection{Sample vs.~Population Eigenvectors}
\label{sec:eigenvectors}

As will be made more apparent in Equation (\ref{eq:phi_link}) below, it is possible to quantify the asymptotic behavior of sample eigenvectors in the general case by selecting a function $g$ of the form $\bone_{(-\infty,\tau)}$ in Equation (\ref{eq:theta}). Let us briefly explain why. \\
First of all, note that each sample eigenvector $u_i$ lies in a space whose dimension is growing towards infinity. Therefore, the only way to know ``where'' it lies is to project it onto a known orthonormal basis that will serve as a reference grid. Given the nature of the problem, the most meaningful choice for this reference grid is the orthonormal basis formed by the population eigenvectors $(v_1,\ldots,v_N)$. Thus we are faced with the task of characterizing the asymptotic behavior of $u_i^*v_j$ for all $i,j=1,\ldots,N$, i.e.~the projection of the sample eigenvectors onto the population eigenvectors. Yet as every eigenvector is identified up to multiplication by a scalar of modulus one, the argument (angle) of $u_i^*v_j$ is devoid of mathematical relevance. Therefore, we can focus instead on its square modulus $\left|u_i^*v_j\right|^2$ without loss of information. \\
Another issue that arises is that of scaling. Indeed as 
$$\frac{1}{N^2}\sum_{i=1}^N\sum_{j=1}^N\left|u_i^*v_j^{}\right|^2
=\frac{1}{N^2}\sum_{i=1}^Nu_i^*\left(\sum_{j=1}^Nv_j^{}v_j^*\right)u_i^{}
=\frac{1}{N^2}\sum_{i=1}^Nu_i^*u_i^{}=\frac{1}{N},$$ we study $N\left|u_i^*v_j\right|^2$ instead, so that its limit does not vanish under large-$N$ asymptotics.\\
The indexing of the eigenvectors also demands special attention as the dimension goes to infinity. We choose to use an indexation system where ``eigenvalues serve as labels for eigenvectors'', that is $u_i$ is the eigenvector associated to the $i^{\rm th}$ largest eigenvalue $\lambda_i$.

\paragraph{}All these considerations lead us to introduce the following key object:
\be \label{defPhi_N}\forall \lambda,\tau\in\mathbb{R},\qquad\Phi_N(\lambda,\tau)
=\frac{1}{N}\sum_{i=1}^N\sum_{j=1}^N\,|u_i^*v_j^{}|^2\;\bone_{[\lambda_i,+\infty)}(\lambda)
\times\bone_{[\tau_j,+\infty)}(\tau).\ee
This bivariate function is right continuous with left-hand limits and nondecreasing in each of its arguments. It also verifies $\lim_{\genfrac{}{}{0pt}{}{\lambda\to-\infty}{\tau\to-\infty}}\Phi_N(\lambda,\tau)=0$ 
and 
$\lim_{\genfrac{}{}{0pt}{}{\lambda\to+\infty}{\tau\to+\infty}}\Phi_N(\lambda,\tau)=1$. 
Therefore, it satisfies the properties of a bivariate cumulative distribution function. 
\begin{rem} Our function
$\Phi_N$ can be compared with the object introduced in \cite{BMP07}:
$\forall \lambda\in\mathbb{R},$ $F^{S_N}_1(\lambda)
=\sum_{i=1}^N|u_i^*x_N^{}|^2\;\bone_{[\lambda_i,+\infty)}(\lambda),$
where $(x_N)_{N=1,2,\ldots}$ is a sequence of nonrandom unit vectors satisfying the non-trivial condition
$x_N^*\left(\Sigma_N-zI\right)^{-1}x_N\to m_H(z).$
This condition is specified so that projecting the sample eigenvectors onto $x_N$ effectively wipes out any signature of non-rotation-invariant behavior. The main difference is that $\Phi_N$ projects the sample eigenvectors onto the population eigenvectors instead.
\end{rem}

From $\Phi_N$ we can extract precise information about the sample eigenvectors. The average of the quantities of interest $N\left|u_i^*v_j\right|^2$ over the sample (resp.~po\-pulation) eigenvectors associated with the sample (resp.~po\-pulation) eigenvalues lying in the interval $[\underline{\lambda},\overline{\lambda}]$ (resp.~$[\underline{\tau},\overline{\tau}]$) is equal to:
\begin{eqnarray}
\lefteqn{\frac{\sum_{i=1}^N\sum_{j=1}^NN\left|u_i^*v_j\right|^2
\bone_{[\underline{\lambda},\overline{\lambda}]}(\lambda_i)
\times\bone_{[\underline{\tau},\overline{\tau}]}(\tau_j)}{\sum_{i=1}^N\sum_{j=1}^N
\bone_{[\underline{\lambda},\overline{\lambda}]}(\lambda_i)
\times\bone_{[\underline{\tau},\overline{\tau}]}(\tau_j)}}\hspace{3cm}\nonumber\\
&=&
\frac{\Phi_N(\overline{\lambda},\overline{\tau})
-\Phi_N(\overline{\lambda},\underline{\tau})
-\Phi_N(\underline{\lambda},\overline{\tau})
+\Phi_N(\underline{\lambda},\underline{\tau})}{
[F_N(\overline{\lambda})-F_N(\underline{\lambda})]\times
[H_N(\overline{\tau})-H_N(\underline{\tau})]},\label{eq:average}
\end{eqnarray}
whenever the denominator is strictly positive. Since $\underline{\lambda}$ and $\overline{\lambda}$ (resp.~$\underline{\tau}$ and $\overline{\tau}$) can be chosen arbitrarily close to each other (as long as the average in Equation (\ref{eq:average}) exists), our goal of characterizing the behavior of sample eigenvectors would be achieved in principle by determining the asymptotic behavior of $\Phi_N$. This can be deduced from Theorem \ref{theo:theta} thanks to the inversion formula for the Stieltjes transform: for all $(\lambda,\tau)\in\R^2$ such that $\Phi_N$ is continuous at $(\lambda,\tau)$
\be
\label{eq:phi_link}
\Phi_N(\lambda,\tau)
=\lim_{\eta\to0^+}\frac{1}{\pi}\int_{-\infty}^\lambda\im\left[\Theta_N^g(\xi+i\eta)\right]d\xi,
\ee
which holds in the special case where $g=\bone_{(-\infty,\tau)}$. We are now ready to state our second main result.
\bt
\label{theo:kernel}
Assume that conditions $(H_1)-(H_4)$ hold true and let $\Phi_N(\lambda,\tau)$ be defined by (\ref{defPhi_N}).
Then there exists a nonrandom bivariate function $\Phi$ such that $\Phi_N(\lambda,\tau)\stackrel{a.s.}{\longrightarrow}\Phi(\lambda,\tau)$ at all points of continuity of $\Phi$. Furthermore, when $\gamma\neq1$, the function $\Phi$ can be expressed as: 
$\forall(\lambda,\tau)\in\mathbb{R}^2,\quad\Phi(\lambda,\tau)
=\int_{-\infty}^\lambda\int_{-\infty}^\tau\varphi(l,t)\,dH(t)\,dF(l),$ where 
\be
\label{eq:kernel}
\forall(l,t)\in\mathbb{R}^2\quad\varphi(l,t)=\left\{
\begin{array}{cl}
\displaystyle\frac{\gamma^{-1}lt}{\left(at-l\right)^2+b^2t^2}
&\mbox{if $l>0$}\\
\displaystyle\frac{1}{(1-\gamma)[1+\breve{m}_{\underline{F}}(0)\,t]}
&\mbox{if $l=0$ and $\gamma<1$}\\
0&\mbox{otherwise,}
\end{array}\right.
\ee
and $a$ (resp.~$b$) is the real (resp.~imaginary) part of $1-\gamma^{-1}-\gamma^{-1}l\,\breve{m}_F(l)$.
\et
Equation (\ref{eq:kernel}) quantifies how the eigenvectors of the sample covariance matrix deviate from those of the population covariance matrix under large-dimensional asymptotics. The result is explicit as a function of $m_F$. 

\begin{figure}[htb]
\centering
\includegraphics{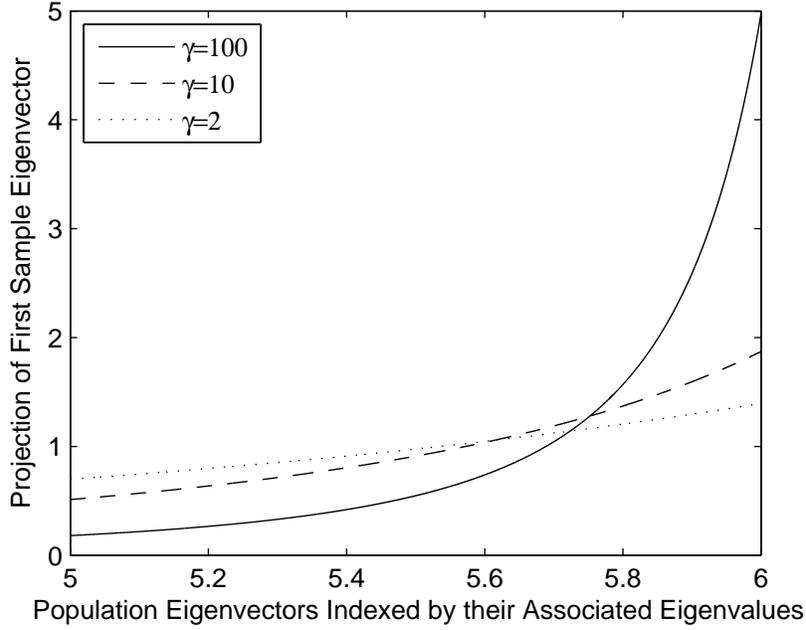}
\caption{\protect\label{fig:eigenvector} Projection of first sample eigenvector onto population eigenvectors (indexed by their associated eigenvalues). We have taken $H'=\bone_{[5,6]}$.}
\end{figure}

To illustrate Theorem \ref{theo:kernel}, we can pick any eigenvector of our choosing, for example the one that corresponds to the first (i.e.~largest) eigenvalue, and plot how it projects onto the population eigenvectors (indexed by their corresponding eigenvalues). The resulting graph is shown in Figure \ref{fig:eigenvector}.
This is a plot of $\varphi(l,t)$ as a function of $t$, for fixed $l$ equal to the supremum of $\supp(F)$. It is the asymptotic equivalent to plotting $N|u_{1}^*v_j|^2$ as a function of $\tau_j$. It looks like a density because, by construction, it must integrate to one. As soon as the sample size starts to drop below $10$ times the number of variables, we can see that the first sample eigenvector starts deviating quite strongly from the first population eigenvectors. This should have precautionary implications for Principal Component Analysis (PCA), where the number of variables is often so large that it is difficult to make the sample size more than ten times bigger. 

Obviously, Equation (\ref{eq:kernel}) would enable us to draw a similar graph for any sample eigenvector (not just the first one), and for any $\gamma$ and $H$ verifying the assumptions of Theorem \ref{theo:kernel}. Preliminary investigations reveal some unexpected patterns. For example: one might have thought that the sample eigenvector associated with the median sample eigenvalue would be closest to the population eigenvector associated with the median population eigenvalue; but in general this is not true.

\subsection{Asymptotically Optimal Bias Correction for the Sample Eigenvalues}
\label{sec:shrinkage}

\begin{figure}[htb]
\centering
\includegraphics{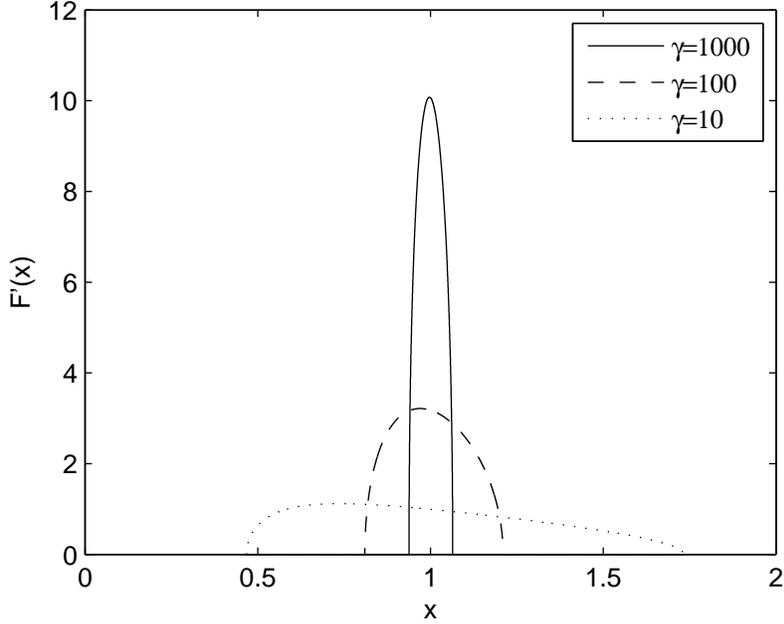}
\caption{\protect\label{fig:dispersion} Limiting density of sample eigenvalues, in the particular case where all the eigenvalues of the population covariance matrix are equal to one. The graph shows excess dispersion of the sample eigenvalues. The formula for this plot comes from solving the Mar\v{c}enko-Pastur equation for $H=\bone_{[1,+\infty)}$.}
\end{figure}

We now bring the two preceding results together to quantify the relationship between the sample covariance matrix and the population covariance matrix \emph{as a whole}. As will be made clear in Equation (\ref{eq:delta_link}) below, this is achieved by selecting the function $g(\tau)=\tau$ in Equation (\ref{eq:theta}). The objective is to see how the sample covariance matrix deviates from the population covariance matrix, and how we can modify it to bring it closer to the population covariance matrix. The main problem with the sample covariance matrix is that its eigenvalues are too dispersed: the smallest ones are biased downwards, and the largest ones upwards. This is most easily visualized when the population covariance matrix is the identity, in which case the limiting spectral e.s.d.~$F$ is known in closed form (see Figure \ref{fig:dispersion}). 
We can see that the smallest and the largest sample eigenvalues are biased away from one, and that the bias decreases in $\gamma$. Therefore, a key concern in multivariate statistics is to find the asymptotically optimal bias correction for the eigenvalues of the sample covariance matrix. As this correction will tend to reduce the dispersion of the eigenvalues, it is often called a \emph{shrinkage} formula.

\begin{figure}[htb]
\centering
\includegraphics{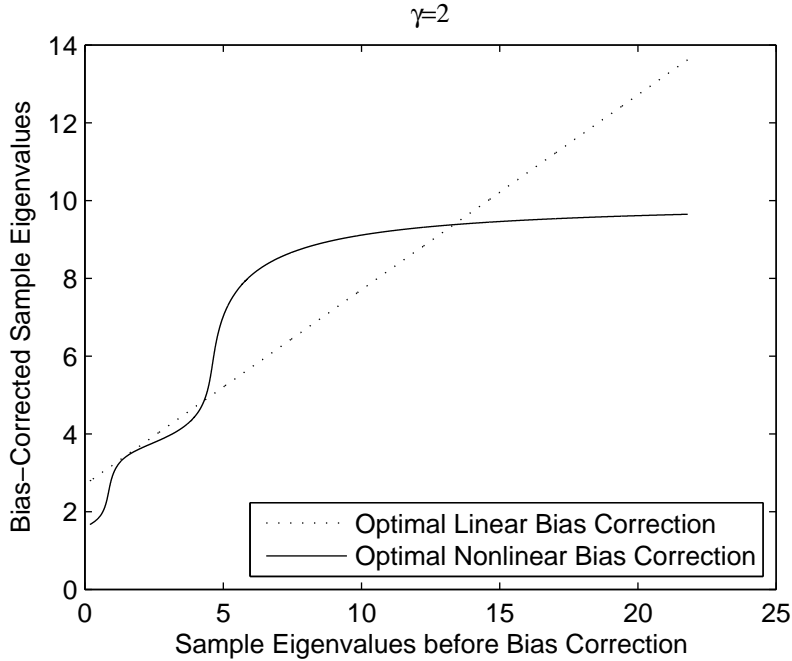}
\caption{\protect\label{fig:lin_nonlin} Comparison of the Optimal Linear vs.~Nonlinear Bias Correction Formul\ae. In this example, the distribution of population eigenvalues $H$ places $20\%$ mass at $1$, $40\%$ mass at $3$ and $40\%$ mass at $10$. The solid line plots $\delta(\lambda)$ as a function of $\lambda$.}
\end{figure}

Ledoit and Wolf \cite{LW04} made some progress along this direction by finding the optimal \emph{linear} shrinkage formula for the sample eigenvalues (projecting $\Sigma_N$ on the two-dimensional subspace spanned by $S_N$ and $I$). However, shrinking the eigenvalues is a highly nonlinear problem (as Figure \ref{fig:lin_nonlin} below will illustrate). Therefore, there is strong reason to believe that finding the optimal \emph{nonlinear} shrinkage formula for the sample eigenvalues would lead to a covariance matrix estimator that further improves upon the Ledoit-Wolf estimator. Theorem \ref{theo:theta} paves the way for such a development.\\
To see how, let us think of the problem of estimating $\Sigma_N$ in general terms. In order to construct an estimator of $\Sigma_N$, we must in turn consider what the eigenvectors and the eigenvalues of this estimator should be. 
Let us consider the eigenvectors first. In the general case where we have no prior information about the orientation of the population eigenvectors, it is reasonable to require that the estimation procedure be invariant with respect to rotation by any $p$-dimensional orthogonal matrix $W$. If we rotate the variables by $W$, then we would ask our estimator to also rotate by the same orthogonal matrix $W$. The class of orthogonally invariant estimators of the covariance matrix is constituted of all the estimators that have the same eigenvectors as the sample covariance matrix (see \cite{P07}, Lemma 5.3). Every rotation-invariant estimator of $\Sigma_N$ is thus of the form: 
$$U_N^{}D_N^{}U_N^*,\quad\mbox{where}\quad D_N=\diag(d_1,\ldots,d_N)\;\mbox{is diagonal,}$$
and where $U_N$ is the matrix whose $\ith$ column is the sample eigenvector $u_i$. This is the class that we consider.\\
Our objective is to find the matrix in this class that is closest to the population covariance matrix. In order to measure distance, we choose the Frobenius norm, defined as:
$\|A\|_F=\sqrt{\tr\left(AA^*\right)}$
for any matrix $A$. Thus we end up with the following optimization problem:
$\min_{D_N\;\mbox{diagonal}}\|U_N^{}D_N^{}U_N^*-\Sigma_N\|_F$.
Elementary matrix algebra shows that its solution is:
$$\widetilde{D}_N=\diag(\widetilde{d}_1,\ldots,\widetilde{d}_N)\quad\mbox{where}\quad\forall i=1,\ldots,N\quad
\widetilde{d}_i=u_i^*\,\Sigma_N\,u_i^{}.$$
The interpretation of $\widetilde{d}_i$ is that it captures how the $\ith$ sample eigenvector $u_i$ relates to the population covariance matrix $\Sigma_N$ as a whole.\\
While $U_N\widetilde{D}_NU_N^*$ does not constitute a \emph{bona fide} estimator (because it depends on the unobservable $\Sigma_N$), new estimators that seek to improve upon the existing ones will need to get as close to $U_N\widetilde{D}_NU_N^*$ as possible. This is exactly the path that led Ledoit and Wolf \cite{LW04} to their improved covariance matrix estimator. Therefore, it is important, in the interest of developing a new and improved estimator, to characterize the asymptotic behavior of $\widetilde{d}_i$ $(i=1,\ldots,N)$. The key object that will enable us to achieve this goal is the nondecreasing function defined by:
\be \label{defDelta_N}\forall x\in\mathbb{R},\quad\Delta_N(x)=\frac{1}{N}\sum_{i=1}^N\widetilde{d}_i\;\bone_{[\lambda_i,+\infty)}(x)
=\frac{1}{N}\sum_{i=1}^Nu_i^*\Sigma_Nu_i\times\bone_{[\lambda_i,+\infty)}(x).
\ee
When all the sample eigenvalues are distinct, it is straightforward to recover the $\widetilde{d}_i$'s from $\Delta_N$:
\be
\label{eq:delta_inverse}
\forall i=1,\ldots,N\qquad\widetilde{d}_i
=\lim_{\varepsilon\to0^+}\frac{\Delta_N(\lambda_i+\varepsilon)-\Delta_N(\lambda_i-\varepsilon)}{
F_N(\lambda_i+\varepsilon)-F_N(\lambda_i-\varepsilon)}.
\ee
The asymptotic behavior of $\Delta_N$ can be deduced from Theorem \ref{theo:theta} in the special case where $g(\tau)=\tau$: for all $x\in\R$ such that $\Delta_N$ continuous at $x$
\be
\label{eq:delta_link}
\Delta_N(x)
=\lim_{\eta\to0^+}\frac{1}{\pi}\int_{-\infty}^x\im\left[\Theta_N^g(\xi+i\eta)\right]d\xi, \quad g(x)\equiv x.
\ee
We are now ready to state our third main result.
\bt
\label{theo:Delta}
Assume that conditions $(H_1)-(H_4)$ hold true and let $\Delta_N$ be defined as in (\ref{defDelta_N}).
There exists a nonrandom function $\Delta$ defined over $\mathbb{R}$ such that $\Delta_N(x)$ converges a.s.~to $\Delta(x)$ for all $x\in\mathbb{R}-\{0\}$. If in addition $\gamma\neq1$, then $\Delta$ can be expressed as: $\forall x\in\mathbb{R},\quad\Delta(x)=\int_{-\infty}^x\delta(\lambda)\,dF(\lambda)$, where
\be \label{def: delta}\forall\lambda\in\mathbb{R},\qquad\delta(\lambda)=\left\{
\begin{array}{cl}
\displaystyle
\frac{\lambda}{\left|1-\gamma^{-1}-\gamma^{-1}\lambda\;\breve{m}_F(\lambda)\right|^2\,}
&\mbox{if $\lambda>0$}\\
\displaystyle\frac{\,\gamma}{(1-\gamma)\,\breve{m}_{\underline{F}}(0)}
&\mbox{if $\lambda=0$ and $\gamma<1$}\\
0&\mbox{otherwise.}
\end{array}\right.\ee
\et
By Equation (\ref{eq:delta_inverse}) the asymptotic quantity that corresponds to $\widetilde{d}_i=u_i^*\Sigma_Nu_i^{}$ is $\delta(\lambda)$, provided that $\lambda$ corresponds to $\lambda_i$. Therefore, the way to get closest to the population covariance matrix (according to the Frobenius norm) would be to divide each sample eigenvalue $\lambda_i$ by the correction factor $|1-\gamma^{-1}-\gamma^{-1}\lambda\,\breve{m}_{F}(\lambda_i)|^2$. This is what we call the optimal nonlinear shrinkage formula or asymptotically optimal bias correction.\footnote{This approach cannot possibly generate a consistent estimator of the population covariance matrix according to the Frobenius norm when $\gamma$ is finite. At best, it could generate a consistent estimator of the {\em projection} of the population covariance matrix onto the space of matrices that have the same eigenvectors as the sample covariance matrix.} Figure \ref{fig:lin_nonlin} shows how much it differs from Ledoit and Wolf's \cite{LW04} optimal linear shrinkage formula. 
In addition, when $\gamma<1$, the sample eigenvalues equal to zero need to be replaced by $\delta(0)=\gamma/[(1-\gamma)\,\breve{m}_{\underline{F}}(0)]$.\\ 
In a statistical context of estimation, $\breve{m}_{F}(\lambda_i)$ and $\breve{m}_{\underline{F}}(0)$ are not known, so they need to be replaced by $\breve{m}_{\widehat{F}}(\lambda_i)$ and $\breve{m}_{\widehat{\underline{F}}}(0)$ respectively, where $\widehat{F}$ is some estimator of the limiting p.d.f.~of sample eigenvalues. Research is currently underway to prove that a covariance matrix estimator constructed in this manner has desirable properties under large-dimensional asymptotics. 

A recent paper \cite{EK08} introduced an algorithm for deducing the population eigenvalues from the sample eigenvalues using the Mar\v{c}enko-Pastur equation. But our objective is quite different, as it is not the population eigenvalues $\tau_i=v_i^*\,\Sigma_N\,v_i^{}$ that we seek, but instead the quantities $\widetilde{d}_i=u_i^*\,\Sigma_N\,u_i^{}$, which represent the diagonal entries of the orthogonal projection (according to the Frobenius norm) of the population covariance matrix onto the space of matrices that have the same eigenvectors as the sample covariance matrix. Therefore the algorithm in \cite{EK08} is better suited for estimating the population eigenvalues themselves, whereas our approach is better suited for estimating the population covariance matrix as a whole.

\begin{figure}[htb]
\centering
\includegraphics{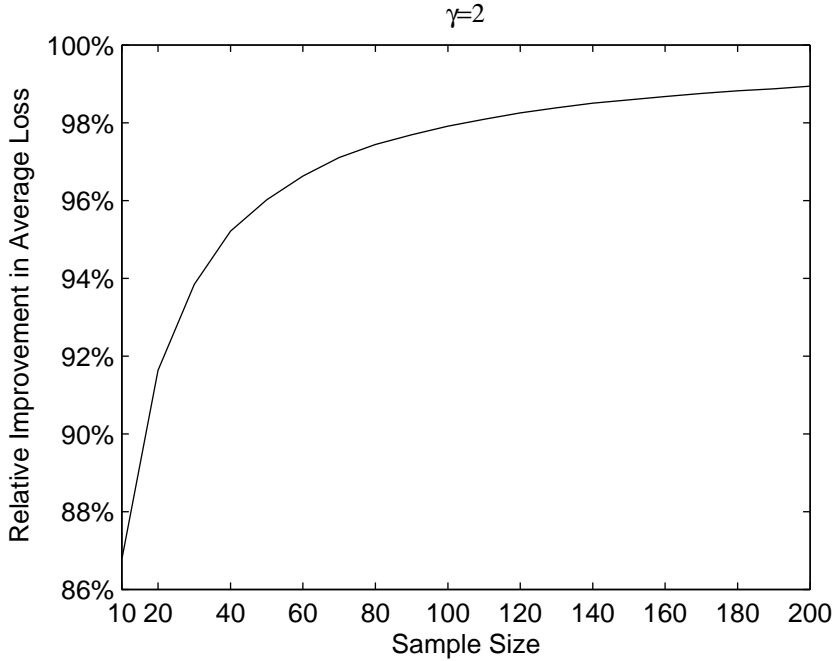}
\caption{\protect\label{fig:convergence} Percentage Relative Improvement in Average Loss (PRIAL) from applying the optimal nonlinear shrinkage formula to the sample eigenvalues. The solid line shows the PRIAL obtained by dividing the $\ith$ sample eigenvalue by the correction factor $|1-\gamma^{-1}-\gamma^{-1}\lambda_i\,\breve{m}_{F}(\lambda_i)|^2$, as a function of sample size. The dotted line shows the PRIAL of the Ledoit-Wolf \protect\cite{LW04} linear shrinkage estimator. For each sample size we generated 10,000 Monte-Carlo simulations using the multivariate Gaussian distribution.  Like in Figure \ref{fig:lin_nonlin}, we used $\gamma=2$ and the distribution of population eigenvalues $H$ placing $20\%$ mass at $1$, $40\%$ mass at $3$ and $40\%$ mass at $10$.}
\end{figure}

Monte-Carlo simulations indicate that applying this bias correction is highly beneficial, even in small samples. We ran 10,000 simulations based on the distribution of population eigenvalues $H$ that places $20\%$ mass at $1$, $40\%$ mass at $3$ and $40\%$ mass at $10$. We kept $\gamma$ constant at 2 while increasing the number of variables from 5 to 100. For each set of simulations, we computed the Percentage Relative Improvement in Average Loss (PRIAL). 
The PRIAL of an estimator $M$ of $\Sigma_N$ is defined as  
$$PRIAL(M)=100 \times \left [1-\frac{\E\left\|M-U_N\widetilde{D}_NU_N^*\right\|_F^2}{\E\left\|S_N-U_N\widetilde{D}_NU_N^*\right\|_F^2}\right ].$$
By construction, the PRIAL of the sample covariance matrix $S_N$ (resp. of $U_N\widetilde{D}_NU_N^*$) is $0\%$ (resp. $100\%$), meaning no improvement (resp. meaning maximum attainable improvement).
For each of the 10,000 Monte-Carlo simulations, we consider $\widetilde{S}_N$, which is the matrix obtained from the sample covariance matrix by keeping its eigenvectors and dividing its $\ith$ eigenvalue by the correction factor $|1-\gamma^{-1}-\gamma^{-1}\lambda_i\,\breve{m}_{F}(\lambda_i)|^2$. The expected 
loss $\E\left\|\widetilde{S}_N-U_N\widetilde{D}_NU_N^*\right\|_F^2$ is estimated by computing its average across the 10,000 Monte-Carlo simulations. Figure \ref{fig:convergence} plots the PRIAL obtained in this way, that is by applying the optimal nonlinear shrinkage formula to the sample eigenvalues. We can see that, even with a modest sample size like $p=40$, we already get $95\%$ of the maximum possible improvement. 

A similar formula can be obtained for the purpose of estimating the \emph{inverse} of the population covariance matrix. To this aim, we set $g(\tau)=1/\tau$ in Equation (\ref{eq:theta}) and define
$$\Psi_N(x):=N^{-1}\sum_{i=1}^Nu_i^*\Sigma_N^{-1}u_i\times\bone_{[\lambda_i,+\infty)}(x), \, \forall x \in \mathbb{R}.$$
\bt\label{theo4b}
Assume that conditions $(H_1)-(H_4)$ are satisfied.
There exists a nonrandom function $\Psi$ defined over $\mathbb{R}$, such that $\Psi_N(x)$ converges a.s.~to $\Psi(x)$ for all $x\in\mathbb{R}-\{0\}$. If in addition $\gamma\neq1$, then $\Psi$ can be expressed as: $\forall x\in\mathbb{R},\quad\Psi(x)=\int_{-\infty}^x\psi(\lambda)\,dF(\lambda)$, where
\be \label{def: psi}\forall\lambda\in\mathbb{R}\quad\psi(\lambda)=\left\{
\begin{array}{cl}
\displaystyle
\frac{1-\gamma^{-1}-2\gamma^{-1}\lambda\,\re\left[\breve{m}_F(\lambda)\right]}{\lambda}
&\mbox{if $\lambda>0$}\\
\displaystyle
\frac{1}{1-\gamma}\,\breve{m}_H(0)-\breve{m}_{\underline{F}}(0)
&\mbox{if $\lambda=0$ and $\gamma<1$}\\
0\;\qquad&\mbox{otherwise.}
\end{array}\right.\ee
\et
Therefore, the way to get closest to the inverse of the population covariance matrix (according to the Frobenius norm) would be to multiply the inverse of each sample eigenvalue $\lambda_i^{-1}$ by the correction factor $1-\gamma^{-1}-2\gamma^{-1}\lambda_i\,\re[\breve{m}_{F}(\lambda_i)]$.  This represents the optimal nonlinear shrinkage formula (or asymptotically optimal bias correction) for the purpose of estimating the inverse covariance matrix. 
Again, in a statistical context of estimation, the unknown $\breve{m}_{F}(\lambda_i)$ needs to be replaced by $\breve{m}_{\widehat{F}}(\lambda_i)$, where $\widehat{F}$ is some estimator of the limiting p.d.f.~of sample eigenvalues. This question is investigated in some work under progress.

\paragraph{}The rest of the paper is organized as follows. Section \ref{sec:prooftheta} contains the proof of Theorem \ref{theo:theta}. Section \ref{sec:proofkernel} contains the proof of Theorem \ref{theo:kernel}. Section \ref{sec:remproofs} is devoted to the proofs of Theorems \ref{theo:Delta} and \ref{theo4b}.

\section{Proof of Theorem \protect\ref{theo:theta}\label{sec:prooftheta}}
The proof of Theorem \protect\ref{theo:theta} follows from an extension of the usual proof of the Mar\v{c}enko-Pastur theorem (see e.g. \cite{Silverstein} and \cite{SilversteinBai}). The latter is based on the Stieltjes transform and, essentially, on a recursion formula. First, we slightly modify this proof to consider more general functionals $\Theta^g_N$ for some polynomial functions $g$. Then we use a standard approximation scheme to extend Theorem \protect\ref{theo:theta} to more general functions $g$. 

First we need to adapt a Lemma from Bai and Silverstein \cite{SilversteinBai}.
\bl
\label{lem: convtr}
Let $Y=(y_1,\ldots,y_N)$ be a random vector with i.i.d.~entries satisfying:
$$\E y_1=0,\quad\E|y_1|^2=1,\quad\E|y_1|^{12}\leq B,$$
where the constant $B$ does not depend on $N$.  Let also $A$ be a given $N\times N$ matrix. Then there exists a constant $K>0$ independent of $N$, $A$ and $Y$ such that:
$$\E\left|YAY^*-\tr(A)\right|^6\leq K\|A\|^6N^3.$$
\el
\paragraph{Proof of Lemma \ref{lem: convtr}}

The proof of Lemma \ref{lem: convtr} directly follows from that of Lemma 3.1 in \cite{SilversteinBai}. Therein the assumption that $\E|y_1|^{12}\leq B$ is replaced with the assumption that $|y_1|\leq \ln N.$ One can easily check that all their arguments carry through if one assumes that the twelfth moment of $y_1$ is uniformly bounded. $\square$

\paragraph{}
Next, we need to introduce some notation. We set $R_N(z)=(S_N-zI)^{-1}$ and define $\Theta_N^{(k)}(z)=N^{-1}\tr[R_N(z)\Sigma^k]$ for all $z\in\C^+$ and integer $k$. Thus, $\Theta_N^{(k)}=\Theta_N^g$ if we take $g(\tau)=\tau^k,\;\;\forall \tau\in\R$. In particular, $\Theta_N^{(0)}=m_{F_N}$. To avoid confusion, the dependency of most of the variables on $N$ will occasionally be dropped from the notation. All convergence statements will be as $N\to\infty$. Conditions $(H_1)-(H_4)$ are assumed to hold throughout.
\bl
\label{lem: convHi}
One has that $\forall z\in\C^+,\quad\Theta_N^{(1)}(z)\stackrel{\rm a.s.}{\longrightarrow}\Theta^{(1)}(z)$ where:
$$\Theta^{(1)}(z)=\frac{\gamma^2}{\gamma-1-z\,m_F(z)}-\gamma.$$
\el
\paragraph{Proof of Lemma \ref{lem: convHi}}
In the first part of the proof, we show that 
$$1+zm_{F_N}(z)=\frac{p}{N}-\frac{1}{N}\sum_{k=1}^p \frac{1}{1+(N/p) \Theta_N^{(1)}(z)}+o(1).$$
Using the a.s.~convergence of the Stieltjes transform $m_{F_N}(z)$, it is then easy to deduce the equation satisfied by $\Theta^{(1)}$ in Lemma \ref{lem: convHi}. Our proof closely follows some of the ideas of \cite{Silverstein} and \cite{SilversteinBai}. Therein the convergence of the Stieltjes transform $m_{F_N}(z)$ is investigated. 

Let us define $C_k=p^{-1/2}\sqrt{\Sigma}\,X_k,$ where $X_k$ is the $\kth$ column of $X$.
Then $S_N =\sum_{k=1}^p C_kC_k^*.$ Using the identity $S_N-zI+zI= \sum_{k=1}^p C_kC_k^*$, one deduces that 
\be \label{eq1}\frac{1}{N}\text{Tr}( I+zR_N(z))=\frac{1}{N}\sum_{k=1}^p C_k^* R_N(z)C_k.\ee
Define now for any integer $1\leq k\leq p$ 
$$R_N^{(k)}(z):=(S_N-C_kC_k^*-zI)^{-1}.$$
By the resolvent identity $R_N(z)-R_N^{(k)}(z)=-R_N(z)C_kC_k^*R_N^{(k)}(z),$ we deduce that
$$C_k^*R_N(z)C_k-C_k^*R_N^{(k)}(z)C_k=-C_k^*R_N(z)C_kC_k^*R_N^{(k)}(z)C_k,$$
which finally gives that
$$C_k^*R_N(z)C_k=\frac{1}{1+C_k^*R_N^{(k)}(z)C_k}C_k^*R_N^{(k)}(z)C_k.$$
Plugging the latter formula into (\ref{eq1}), one can write that
\be \label{eq2}
1+zm_{F_N}(z)=\frac{p}{N}-\frac{1}{N}\sum_{k=1}^p \frac{1}{1+C_k^*R_N^{(k)}(z)C_k}.\ee
We will now use the fact that $R_N^{(k)}$ and $C_k$ are independent random matrices to estimate the asymptotic behavior of the last sum in (\ref{eq2}). Using Lemma \ref{lem: convtr}, we deduce that 
\be \label{conCk}\max_{k\in \{1, \ldots, p\}} \Big |C_k^* R_N^{(k)}C_k- \frac{1}{p}\tr\left(R_N^{(k)}\Sigma\right)\Big |\stackrel{\text{a.s.}}{\to}0, \ee
as $N \to \infty.$
Furthermore, using Lemma 2.6 in Silverstein and Bai (1995), one also has that 
\be \label{starstar}\frac{1}{p} \Big|\tr\left[ \left(R_N-R_N^{(k)}\right)\Sigma\right] \Big|\leq \frac{||\Sigma ||}{py}.\ee
Thus using (\ref{starstar}), (\ref{conCk}) and (\ref{eq2}), one can write that 
\be\label{mastereq}
1+zm_{F_N}(z)=\frac{p}{N}-\frac{1}{N}\sum_{k=1}^p \frac{1}{1+(N/p) \Theta_N^{(1)}(z)}+\delta_N,
\ee
where the error term $\delta_N$ is given by $\delta_N=\delta_N^1+\delta_N^2$ with 
$$\delta_N^1=\frac{1}{N}\sum_{k=1}^p \frac{\frac{1}{p}\text{Tr}\left ( (R_N-R_N^{(k)}) \Sigma\right) }{(1+\frac{1}{p}\text{Tr}(R_N\Sigma))(1+\frac{1}{p}\text{Tr}(R_N^{(k)}\Sigma))}$$
and $$\delta_N^2=-\frac{1}{N}\sum_{k=1}^p \frac{C_k^* R_N^{(k)}C_k- \text{Tr}(R_N^{(k)}\Sigma)}{(1+\frac{1}{p}\text{Tr}(R_N^{(k)}\Sigma))(1+\frac{1}{p}C_k^*R_N^{(k)}C_k)} .$$
We will now use (\ref{starstar}) and (\ref{conCk}) to show that $\delta_N$ a.s. converges to $0$ as $N \to \infty.$ It is known that $F_N$ converges a.s.~to the distribution $F$ given by the Mar\v{c}enko-Pastur equation (and has no subsequence vaguely convergent to 0). It is proven in Silverstein and Bai (1995) that under these assumptions, there exists $m>0$ such that $\inf_N F_N([-m,m])>0.$
In particular, there exists $\delta >0$ such that
$$\inf_N \im\left[\int \frac{1}{\lambda-z}dF_N(\lambda)\right]\geq \int \frac{y}{2\lambda^2+2x^2+y^2}dF_{N}(\lambda)\geq \delta.$$
From this, we deduce that 
$$\left|1+\frac{1}{p}\tr(\Sigma R_N)\right|\geq \im\left[\frac{1}{p} \tr(\Sigma R_N)\right]\geq \frac{h_1}{\gamma}\delta.$$
Using (\ref{starstar}) we also get that 
$$\left|1+\frac{1}{p}\tr\left(\Sigma R_N^{(k)}\right)\right|\geq \im \left[\frac{1}{p}\tr\left(\Sigma R_N^{(k)}\right)\right]\geq \frac{h_1}{2\gamma}\delta.$$
We first consider $\delta_N^1$.
Thus one has that 
\be
\left|\delta_N^1\right|\leq \frac{2 ||\Sigma|| \gamma^2}{Ny h_1^2 \delta^2} =O(1/N).
\ee
We now turn to $\delta_N^2$. Using the a.s. convergence (\ref{conCk}), it is not hard to deduce that 
$$\delta_N^2 \to 0, \text{ a.s.}$$
This completes the proof of Lemma \ref{lem: convHi}. $\square$

\bl
\label{lem: convHk}
For every $k=1,2,\ldots$ the limit $\lim_{N \to \infty}\Theta_N^{(k)}(z):= \Theta^{(k)}(z)$ exists and satisfies the recursion equation
\be \label{convHk}
\forall z\in\C^+,\quad
 \Theta^{(k+1)}(z)=\left [z\Theta^{(k)}(z)+\int_{-\infty}^{+\infty} \tau^k dH(\tau) \right ]\times\left [1+\frac{1}{\gamma}\Theta^{(1)}(z)\right].\ee
\el
\paragraph{Proof of Lemma \ref{lem: convHk}}The proof is inductive, so we assume that formula (\ref{convHk}) holds for any integer smaller than or equal to $q$ for some given integer $q$.
We start from the formula 
$$\tr\left (\Sigma^q+z\Sigma^q R_N(z)\right)=\tr\left ( \Sigma^q R_N (z) S_N\right)=\sum_{k=1}^p C_k^*\Sigma^q R_N (z)C_k.$$
Using once more the resolvent identity, one gets that 
$$C_k^* \Sigma^q R_N(z)C_k=\frac{C_k^* \Sigma^q R_N^{(k)}(z) C_k}{1+C_k^* R_N^{(k)}(z) C_k},$$
which yields that 
\be \frac{1}{N}\tr\left (\Sigma^q+z\Sigma^q R_N(z)\right)=\sum_{k=1}^p \frac{C_k^* \Sigma^q R_N^{(k)}(z) C_k}{1+C_k^* R_N^{(k)}(z) C_k}.\ee
It is now an easy consequence of the arguments developed in the case where $q=0$ to check that 
$$\max_{k \in \{1, \ldots, p\}}\Big |C_k^*  R_N^{(k)}(z) C_k- \tr\left (\Sigma R_N(z)\right)\Big |+\Big |C_k^* \Sigma^q R_N^{(k)}(z) C_k- \tr\left (\Sigma^{q+1}R_N(z)\right)\Big |$$
converges a.s.~to zero. Using the recursion assumption that $\lim_{N \to \infty}\Theta_N^{(k)}(z)$ exists, $\forall k \leq q$, one can deduce that 
$\lim_{N \to \infty}\Theta_N^{(q+1)}(z)$ exists and that the limit $\Theta^{(q+1)}(z)$ satisfies 
$$\left [z\Theta^{(q)}(z)+\int_{-\infty}^{+\infty} \tau^q dH(\tau) \right ]\times\left [1+\frac{1}{\gamma}\Theta^{(1)}(z)\right] =\Theta^{(q+1)}(z).$$
This finishes the proof of Lemma \ref{lem: convHk}. $\square$

\bl
\label{lem:polynomial}
Theorem \ref{theo:theta} holds when the function $g$ is a polynomial.
\el
\paragraph{Proof of Lemma \ref{lem:polynomial}}Given the linearity of the problem, it is sufficient to prove that Theorem \ref{theo:theta} holds when the function $g$ is of the form: $\forall\tau\in\R, \quad g(\tau)=\tau^k$, for any nonnegative integer $k$. In the case where $k=0$, this is a direct consequence of Theorem 1.1 in \cite{Silverstein}.

The existence of a function $\Theta^{(k)}$ defined on $\C^+$ such that $\Theta_N^{(k)}(z)\stackrel{\rm a.s.}{\longrightarrow}\Theta^{(k)}(z)$ for all $z\in\C^+$ is established by Lemma \ref{lem: convHi} for $k=1$ and by Lemma \ref{lem: convHk} for $k=2,3,\ldots$ Therefore, all that remains to be shown is that Equation (\ref{eq:theta}) holds for $k=1,2,\ldots$ 

We will first show it for $k=1$. From the original Mar\v{c}enko-Pastur equation we know that:
\be 
1+zm_F(z)=\int_{-\infty}^{+\infty}
\frac{\tau\left[1-\gamma^{-1}-\gamma^{-1}z\,m_F(z)\right]}{\tau\left[1-\gamma^{-1}-\gamma^{-1}z\,m_F(z)\right]-z}dH(\tau).
\label{eq:eq1}
\ee
From Lemma \ref{lem: convHi} we know that:
$$ \Theta^{(1)}(z)=\frac{\gamma^2}{\gamma-1-z\,m_F(z)}-\gamma=\frac{1+z\,m_F(z)}{1-\gamma^{-1}-\gamma^{-1}z\,m_F(z)},$$
yielding that 
\be 1+zm_F(z)=\frac{\Theta^{(1)}(z)}{1+\gamma^{-1}\Theta^{(1)}(z)}.\label{eq:eq2}
\ee
Combining Equations (\ref{eq:eq1}) and (\ref{eq:eq2}) yields:
\be
\label{eq:eq3}
\int_{-\infty}^{+\infty}
\frac{\tau\left[1-\gamma^{-1}-\gamma^{-1}z\,m_F(z)\right]}{\tau\left[1-\gamma^{-1}-\gamma^{-1}z\,m_F(z)\right]-z}\,dH(\tau)
=\frac{\Theta^{(1)}(z)}{1+\gamma^{-1}\Theta^{(1)}(z)}.
\ee
From Lemma \ref{lem: convHi}, we also know that: \be 
1+\gamma^{-1}\Theta^{(1)}(z)=\frac{1}{1-\gamma^{-1}-\gamma^{-1}z\,m_F(z)}.\label{eq:eq4}\ee
Putting together Equations (\ref{eq:eq3}) and (\ref{eq:eq4}) yields the simplification:
$$
\Theta^{(1)}(z)=\int_{-\infty}^{+\infty}
\frac{1}{\tau\left[1-\gamma^{-1}-\gamma^{-1}z\,m_F(z)\right]-z}\,\tau\,dH(\tau),
$$
which establishes that Equation (\ref{eq:theta}) holds when $g(\tau)=\tau,$ $\forall\tau\in\R$.

We now show by induction that Equation (\ref{eq:theta}) holds when $g(\tau)=\tau^k$ for $k=2,3,\ldots$ Assume that we have proven it for $k-1$. Thus the recursion hypothesis is that:
\be
\label{eq:eq5}
\Theta^{(k-1)}(z)=\int_{-\infty}^{+\infty}
\frac{1}{\tau\left[1-\gamma^{-1}-\gamma^{-1}z\,m_F(z)\right]-z}\,\tau^{k-1}\,dH(\tau).
\ee
From Lemma \ref{lem: convHk} we know that:
\be
\label{eq:eq6}
\Theta^{(k)}(z)=\left [z\Theta^{(k-1)}(z)+\int_{-\infty}^{+\infty} \tau^{k-1} dH(\tau) \right ]\times\left [1+\frac{1}{\gamma}\Theta^{(1)}(z)\right].
\ee
Combining Equations (\ref{eq:eq5}) and (\ref{eq:eq6}) yields:
\begin{eqnarray}
\frac{\Theta^{(k)}(z)}{1+\frac{1}{\gamma}\Theta^{(1)}(z)}
&=&z\Theta^{(k-1)}(z)+\int_{-\infty}^{+\infty} \tau^{k-1} dH(\tau)\nonumber\\
&=&
\int_{-\infty}^{+\infty}\left\{
\frac{z}{\tau\left[1-\gamma^{-1}-\gamma^{-1}z\,m_F(z)\right]-z}+1\right\}\tau^{k-1}\,dH(\tau)\nonumber\\
&=&
\int_{-\infty}^{+\infty}
\frac{1-\gamma^{-1}-\gamma^{-1}z\,m_F(z)}{\tau\left[1-\gamma^{-1}-\gamma^{-1}z\,m_F(z)\right]-z}\,\tau^k\,dH(\tau).\label{eq:eq7}
\end{eqnarray}
Putting together Equations (\ref{eq:eq4}) and (\ref{eq:eq7}) yields the simplification:
$$\Theta^{(k)}(z)=\int_{-\infty}^{+\infty}
\frac{1}{\tau\left[1-\gamma^{-1}-\gamma^{-1}z\,m_F(z)\right]-z}\,\tau^k\,dH(\tau),$$
which proves that the desired assertion holds for $k$. Therefore, by induction, it holds for all $k=1,2,3,\ldots$ This completes the proof of Lemma \ref{lem:polynomial}. $\square$

\bl
\label{lem:continuous}
Theorem \ref{theo:theta} holds for any function $g$ that is continuous on $[h_1,h_2]$.
\el
\paragraph{Proof of Lemma \ref{lem:continuous}}We shall deduce this from Lemma \ref{lem:polynomial}. Let $g$ be any function that is continuous on $[h_1,h_2]$. By the Weierstrass approximation theorem, there exists a sequence of polynomials that converges to $g$ uniformly on $[h_1,h_2]$. By Lemma \ref{lem:polynomial}, Theorem \ref{theo:theta} holds for every polynomial in the sequence. Therefore it also holds for the limit $g$. $\square$

\paragraph{}
We are now ready to prove Theorem \ref{theo:theta}. We shall prove it by induction on the number $k$ of points of discontinuity of the function $g$ on the interval $[h_1,h_2]$. The fact that it holds for $k=0$ has been established by Lemma \ref{lem:continuous}. Let us assume that it holds for some $k$. Then consider any bounded function $g$ which has $k+1$ points of discontinuity on $[h_1,h_2]$. Let $\nu$ be one of these $k+1$ points of discontinuity. Construct the function: $\forall x\in[h_1,h_2], \quad \rho(x)=g(x)\times(x-\nu)$. The function $\rho$ has $k$ points of discontinuity on $[h_1,h_2]$: all the ones that $g$ has, except $\nu$. Therefore, by the recursion hypothesis, 
$\Theta_N^\rho(z)=N^{-1}\tr\left[(S_N-zI)^{-1}\rho(\Sigma_N)\right]$ converges a.s.~to 
\be
\label{eq:eq8}
\Theta^\rho(z)
=\int_{-\infty}^{+\infty}\frac{1}{\tau\left[1-\gamma^{-1}-\gamma^{-1}z\,m_F(z)\right]-z}\,\rho(\tau)\,dH(\tau)
\ee
for all $z\in\C^+$. It is easy to adapt the arguments developed in the proof of Lemma \ref{lem: convHk} to show that $\lim_{N\to\infty}\Theta_N^g(z)$ exists (as $g$ is bounded) and is equal to:
\be
\label{eq:eq9}
\Theta^g(z)=\frac{\Theta^\rho(z)-\left[1+\gamma^{-1}\Theta^{(1)}(z)\right]\int_{-\infty}^{+\infty}g(\tau)dH(\tau)}{
z\left[1+\gamma^{-1}\Theta^{(1)}(z)\right]-\nu}
\ee
for all $z\in\C^+$. Plugging Equation (\ref{eq:eq8}) into Equation (\ref{eq:eq9}) yields:
$$
\Theta^g(z)=\frac{
\int_{-\infty}^{+\infty}\left\{
\frac{\tau-\nu}{\tau\left[1-\gamma^{-1}-\gamma^{-1}z\,m_F(z)\right]-z}
-\left[1+\gamma^{-1}\Theta^{(1)}(z)\right]\right\}\,g(\tau)\,dH(\tau)}{
z\left[1+\gamma^{-1}\Theta^{(1)}(z)\right]-\nu}.
$$
Using Equation (\ref{eq:eq4}) we get:
\begin{eqnarray*}
\Theta^g(z)&=&\frac{\int_{-\infty}^{+\infty}\left\{
\frac{\tau-\nu}{\tau\left[1-\gamma^{-1}-\gamma^{-1}z\,m_F(z)\right]-z}
-\frac{1}{1-\gamma^{-1}-\gamma^{-1}z\,m_F(z)}\right\}\,g(\tau)\,dH(\tau)}{
\frac{z}{1-\gamma^{-1}-\gamma^{-1}z\,m_F(z)}-\nu}\\
&=&\frac{\int_{-\infty}^{+\infty}
\frac{z-\nu\left[1-\gamma^{-1}-\gamma^{-1}z\,m_F(z)\right]}{
\left\{\tau\left[1-\gamma^{-1}-\gamma^{-1}z\,m_F(z)\right]-z\right\}\times
\left[1-\gamma^{-1}-\gamma^{-1}z\,m_F(z)\right]}\,g(\tau)\,dH(\tau)}{
\frac{z-\nu\left[1-\gamma^{-1}-\gamma^{-1}z\,m_F(z)\right]}{
1-\gamma^{-1}-\gamma^{-1}z\,m_F(z)}}\\
&=&\int_{-\infty}^{+\infty}\frac{1}{\tau\left[1-\gamma^{-1}-\gamma^{-1}z\,m_F(z)\right]-z}\,g(\tau)\,dH(\tau),
\end{eqnarray*}
which means that Equation (\ref{eq:theta}) holds for $g$. Therefore, by induction, Theorem \ref{theo:theta} holds for any bounded function $g$ with a finite number of discontinuities on $[h_1,h_2]$. $\square$

\section{Proof of Theorem \protect\ref{theo:kernel}\label{sec:proofkernel}}

At this stage, we need to establish two Lemmas that will be of general use for deriving implications from Theorem \ref{theo:theta}.
\bl
\label{lem:omega}
Let $g$ denote a (real-valued) bounded function defined on $[h_1,h_2]$ with finitely many points of discontinuity. Consider the function $\Omega^g_N$ defined by:
$$\forall x\in\R, \,\Omega^g_N(x)
=\frac{1}{N}\sum_{i=1}^N\bone_{[\lambda_i,+\infty)}(x)\sum_{j=1}^N\left|u_i^*v_j\right|^2\times g(\tau_j).$$
Then there exists a nonrandom function $\Omega^g$ defined on $\R$ such that $\Omega^g_N(x)\stackrel{\rm a.s.}{\to}\Omega^g(x)$ at all points of continuity of $\Omega^g$. Furthermore, 
\be
\label{eq:omega}
\Omega^g(x)=\lim_{\eta\to0^+}\frac{1}{\pi}\int_{-\infty}^x\im\left[\Theta^g\left(\lambda+i\eta\right)\right]d\lambda
\ee
for all $x$ where $\Omega^g$ is continuous.
\el 
\paragraph{Proof of Lemma \ref{lem:omega}}
The Stieltjes transform of $\Omega^g_N$ is the function $\Theta^g_N$ defined by Equation (\ref{eq:general}). From Theorem \ref{theo:theta}, we know that there 
exists a nonrandom function $\Theta^g$ defined over $\mathbb{C}^+$ such that $\Theta_N^g(z)\stackrel{\rm a.s.}{\to}\Theta^g(z)$ for all $z\in\mathbb{C}^+$. Therefore, Silverstein and Bai's \cite{SilversteinBai} Equation (2.5) implies that: $\lim_{N\to\infty}\Omega_N(x)\equiv\Omega^g(x)$ exists for all $x$ where $\Omega^g$ is continuous. Furthermore,
the Stieltjes transform of $\Omega^g$ is $\Theta^g$. Then Equation (\ref{eq:omega}) is simply the inversion formula for the Stieltjes transform. $\square$

\bl
\label{lem:exchange}
Under the assumptions of Lemma \ref{lem:omega}, if $\gamma>1$ then for all $(x_1,x_2)\in\R^2$:
\be
\label{eq:exchange}
\Omega^g(x_2)-\Omega^g(x_1)=\frac{1}{\pi}\int_{x_1}^{x_2}\lim_{\eta\to0^+}\im\left[\Theta^g(\lambda+i\eta)\right]d\lambda.
\ee
If $\gamma<1$ then Equation (\ref{eq:exchange}) holds for all $(x_1,x_2)\in\R^2$ such that $x_1x_2>0$.
\el
\paragraph{Proof of Lemma \ref{lem:exchange}} One can first note that
$\lim_{z\in\C^+\to x}\im\left[\Theta^g(z)\right]\equiv\im\left[\Theta^g(x)\right]$ exists for all $x\in\R$ (resp.~all  $x\in\R-\{0\}$) in the case where $\gamma>1$ (resp.~$\gamma<1$). This is obvious if $x\in \supp (F)$. In the case where $x \notin \supp (F)$, then it can be deduced from Theorem 4.1 in \cite{SC95} that 
$\dfrac{x}{1-\gamma^{-1}(1+x \breve m_F (x))}\notin \supp (H)$, which ensures the desired result.
Now $\Theta^g$ is the Stieltjes transform of $\Omega^g$. Therefore, Silverstein and Choi's \cite{SC95} Theorem 2.1 implies that:
$$\mbox{$\Omega^g$ is differentiable at $x$ and its derivative is:}
\;\,\frac{1}{\pi}\,\im\left[\Theta^g(x)\right]$$
for all $x\in\R$ (resp.~all $x\in\R-\{0\}$) in the case where $\gamma>1$ (resp.~$\gamma<1$). When we integrate, we get Equation (\ref{eq:exchange}). $\square$
\paragraph{}
We are now ready to proceed with the proof of Theorem \ref{theo:kernel}. Let $\tau\in\R$ be given and take $g=\bone_{(-\infty,\tau)}$. Then we have:
$$\forall z\in\C^+,\quad\Theta^{\bone_{(-\infty,\tau)}}_N(z)
=\frac{1}{N}\sum_{i=1}^N\frac{1}{\lambda_i-z}\sum_{j=1}^N\left|u_i^*v_j\right|^2\times\bone_{(-\infty,\tau)}.$$
Since the function $g=\bone_{(-\infty,\tau)}$ has a single point of discontinuity (at $\tau$), Theorem \ref{theo:theta} implies that $\forall z\in\C^+,\quad\Theta^{\bone_{(-\infty,\tau)}}_N(z)\stackrel{\rm a.s.}{\to}\Theta^{\bone_{(-\infty,\tau)}}(z)$, where:
\be
\label{eq:eq21}
\forall z\in\C^+, \quad\Theta^{\bone_{(-\infty,\tau)}}(z)
=\int_{-\infty}^{\tau}\frac{1}{t\left[1-\gamma^{-1}-\gamma^{-1}z\,m_F(z)\right]-z}\,dH(t).
\ee
Remember from Equation (\ref{eq:phi_link}) that:
$$\Phi_N(\lambda,\tau)
=\lim_{\eta\to0^+}\frac{1}{\pi}\int_{-\infty}^\lambda\im\left[\Theta_N^{\bone_{(-\infty,\tau)}}(l+i\eta)\right]dl.$$
Therefore, by Lemma \ref{lem:omega}, $\lim_{N\to\infty}\Phi_N(\lambda,\tau)$ exists and is equal to:
\be
\label{eq:eq22}
\Phi(\lambda,\tau)=\lim_{\eta\to0^+}\frac{1}{\pi}\int_{-\infty}^\lambda\im\left[\Theta^{\bone_{(-\infty,\tau)}}(l+i\eta)\right]dl,
\ee
for every $(\lambda,\tau)\in\R^2$ where $\Phi$ is continuous. We first evaluate $\Phi(\lambda,\tau)$ in the case where $\gamma>1$, so that the limiting e.s.d.~$F$ is continuously differentiable on all of $\R$. Plugging (\ref{eq:eq21}) into (\ref{eq:eq22}) yields:
\begin{eqnarray}
\Phi(\lambda,\tau)
&=&\lim_{\eta\to0^+}\frac{1}{\pi}\int_{-\infty}^\lambda\im\left\{\int_{-\infty}^{\tau}\frac{1}{t\left[a(l,\eta)+ib(l,\eta)\right]-l-i\eta}\,dH(t)\right\}dl\nonumber\\
&=&\frac{1}{\pi}\int_{-\infty}^\lambda\int_{-\infty}^{\tau}\lim_{\eta\to0^+}\im\left\{\frac{1}{t\left[a(l,\eta)+ib(l,\eta)\right]-l-i\eta}\right\}dH(t)\,dl,\label{eq:eq23}
\end{eqnarray}
where $a(l,\eta)+ib(l,\eta)=1-\gamma^{-1}-\gamma^{-1}(l+i\eta)\,m_F(l+i\eta)$. The last equality follows from Lemma \ref{lem:exchange}. Notice that:
$$\im\left\{\frac{1}{t\left[a(l,\eta)+ib(l,\eta)\right]-l-i\eta}\right\}
=\frac{\eta-b(l,\eta)t}{\left[a(l,\eta)t-l\right]^2+\left[b(l,\eta)t-\eta\right]^2}.$$
Taking the limit as $\eta\to0^+$, we get:
$$
a(l,\eta)\longrightarrow a=\re\left[1-\frac{1}{\gamma}-\frac{l\breve{m}_F(l)}{\gamma}\right],\;
b(l,\eta)\longrightarrow b=\im\left[1-\frac{1}{\gamma}-\frac{l\breve{m}_F(l)}{\gamma}\right].$$
The inversion formula for the Stieltjes transform implies: $\forall l\in\R,$  $F'(l)=\frac{1}{\pi}\im\left[\breve{m}_F(l)\right]$, therefore $b=-\pi\gamma^{-1}lF'(l)$. Thus we have:
\be
\label{eq:eq24}
\lim_{\eta\to0^+}\im\left\{\frac{1}{t\left[a(l,\eta)+ib(l,\eta)\right]-l-i\eta}\right\}
=\frac{\pi\gamma^{-1}lt}{(at-l)^2+b^2t^2}\times F'(l).
\ee
Plugging Equation (\ref{eq:eq24}) back into Equation (\ref{eq:eq23}) yields that:
$$\Phi(\lambda,\tau)
=\int_{-\infty}^\lambda\int_{-\infty}^{\tau}\frac{\gamma^{-1}lt}{(at-l)^2+b^2t^2}\,dH(t)\,dF(l),
$$
which was to be proven. This completes the proof of Theorem \ref{theo:kernel} in the case where $\gamma>1$.

In the case where $\gamma<1$, much of the arguments remain the same, except for an added degree of complexity due to the fact that the limiting e.s.d.~$F$ has a discontinuity of size $1-\gamma$ at zero. This is handled by using the following three Lemmas.
\bl
\label{lem:neighborhood}
If $\gamma\neq 1$, $F$ is constant over the interval $\left(0,(1-\frac{1}{\sqrt \gamma})^2 h_1\right)$.
\el
\paragraph{Proof of Lemma \ref{lem:neighborhood}}If $H$ placed all its weight on $h_1$, then we could solve the Mar\v{c}enko-Pastur equation explicitly for $F$, and the infimum of the support of the limiting e.s.d.~of nonzero sample eigenvalues would be equal to $(1-\gamma^{-1/2})^2\times h_1$. Since, by Assumption $(H_4)$, $H$ places all its weight on points greater than or equal to $h_1$, the infimum of the support of the limiting e.s.d.~of nonzero sample eigenvalues has to be greater than or equal to $(1-\gamma^{-1/2})^2\times h_1$ (see Equation (1.9b) in Bai and Silverstein \cite{BS04}). Therefore, $F$ is constant over the open interval $\left(0,(1-\gamma^{-1/2})^2\times h_1\right)$. $\square$

\bl
\label{lem:zero} Let $\kappa>0$ be a given real number.
Let $\mu$ be a complex holomorphic function defined on the set $\{z\in\C^+:\re[z]\in(-\kappa,\kappa)\}$. If $\mu(0)\in\R$ then:
$$\lim_{\varepsilon\to0^+}\left\{\lim_{\eta\to0^+}\frac{1}{\pi}\int_{-\varepsilon}^{+\varepsilon}
\im\left[-\frac{\mu(\xi+i\eta)}{\xi+i\eta}\right]d\xi\right\}=\mu(0).$$
\el
\paragraph{Proof of Lemma \ref{lem:zero}}For all $\varepsilon$ in $(0,\kappa),$ we have:
\begin{eqnarray}
\lim_{\eta\to0^+}\frac{1}{\pi}\int_{-\varepsilon}^{+\varepsilon}
\im\left[-\frac{1}{\xi+i\eta}\right]d\xi
&=&\lim_{\eta\to0^+}\frac{1}{\pi}\int_{-\epsilon}^{+\epsilon}
\frac{\eta}{\xi^2+\eta^2}d\xi\nonumber\\
&=&\lim_{\eta\to0^+}\frac{1}{\pi}
\left[\arctan \left (\frac{\varepsilon}{\eta}\right )-\arctan\left (\frac{-\varepsilon}{\eta}\right )\right]\nonumber\\
&=&1.\label{eq:arctan}
\end{eqnarray}
Since $\mu$ is continuously differentiable, there exist $\delta >0$, $\beta >0$ such that $|\mu'(z)|\leq \beta, \forall z, |z|\leq \delta$. Using Taylor's theorem, we get that
$
\left|\mu(z)-\mu(0)\right|\leq\beta|z|, \, \forall |z|\leq \delta.$
Now we can perform the following decomposition:
\begin{eqnarray*}
\lefteqn{\lim_{\varepsilon\to0^+}\left\{\lim_{\eta\to0^+}\frac{1}{\pi}\int_{-\varepsilon}^{+\varepsilon}
\im\left[-\frac{\mu(\xi+i\eta)}{\xi+i\eta}\right]d\xi\right\}}\hspace{2cm}\\
&=&\lim_{\varepsilon\to0^+}\left\{\lim_{\eta\to0^+}\frac{1}{\pi}\int_{-\varepsilon}^{+\varepsilon}
\im\left[-\frac{\mu(\xi+i\eta)-\mu(0)+\mu(0)}{\xi+i\eta}\right]d\xi\right\}\\
&=&\mu(0)\lim_{\varepsilon\to0^+}\left\{\lim_{\eta\to0^+}\frac{1}{\pi}\int_{-\varepsilon}^{+\varepsilon}
\im\left[-\frac{1}{\xi+i\eta}\right]d\xi\right\}\\
&&{}+\lim_{\varepsilon\to0^+}\left\{\lim_{\eta\to0^+}\frac{1}{\pi}\int_{-\varepsilon}^{+\varepsilon}
\im\left[-\frac{\mu(\xi+i\eta)-\mu(0)}{\xi+i\eta}\right]d\xi\right\}\\
&=&\mu(0)+\lim_{\varepsilon\to0^+}\left\{\lim_{\eta\to0^+}\frac{1}{\pi}\int_{-\varepsilon}^{+\varepsilon}
\im\left[-\frac{\mu(\xi+i\eta)-\mu(0)}{\xi+i\eta}\right]d\xi\right\},
\end{eqnarray*}
where the last equality follows from Equation (\ref{eq:arctan}). The second term vanishes because:
\begin{eqnarray*}
\lefteqn{\left|\lim_{\varepsilon\to0^+}\left\{\lim_{\eta\to0^+}\frac{1}{\pi}\int_{-\varepsilon}^{+\varepsilon}
\im\left[-\frac{\mu(\xi+i\eta)-\mu(0)}{\xi+i\eta}\right]d\xi\right\}\right|}\hspace{2cm}\\
&\leq&\lim_{\varepsilon\to0^+}\left\{\lim_{\eta\to0^+}\frac{1}{\pi}\int_{-\varepsilon}^{+\varepsilon}
\left|\frac{\mu(\xi+i\eta)-\mu(0)}{\xi+i\eta}\right|d\xi\right\}\\
&\leq&\lim_{\varepsilon\to0^+}\left\{\lim_{\eta\to0^+}\frac{1}{\pi}\int_{-\varepsilon}^{+\varepsilon}
\beta\,d\xi\right\}=0.
\end{eqnarray*} This yields Lemma \ref{lem:zero}. $\square$ 

\bl
\label{lem:zero_g}
Assume that $\gamma<1$. Let $g$ be a (real-valued) bounded function defined on $[h_1,h_2]$ with finitely many points of discontinuity. Then:
\begin{eqnarray*}
\lefteqn{\lim_{\varepsilon\to0^+}\lim_{\eta\to0^+}\frac{1}{\pi}\int_{-\varepsilon}^{+\varepsilon}
\int_{-\infty}^{+\infty}\im\left\{
\frac{\scriptstyle g(\tau)}{\scriptstyle\tau\left[1-\gamma^{-1}-\gamma^{-1}(\xi+i\eta)\,m_F(\xi+i\eta)\right]-\xi-i\eta}\right\}dH(\tau)d\xi}
\hspace{7cm}\\
&=&\int_{-\infty}^{+\infty}\frac{g(\tau)}{1+\breve{m}_{\underline{F}}(0)\tau}\,dH(\tau),
\end{eqnarray*}
where $\underline{F}=\left(1-\gamma^{-1}\right)\bone_{[0,+\infty)}+\gamma^{-1}F$, and 
$\breve{m}_{\underline{F}}(0)=\lim_{z\in\mathbb{C}^+\to0}m_{\underline{F}}(z)$.
\el
\paragraph{Proof of Lemma \ref{lem:zero_g}} One has that 
\begin{eqnarray}
\forall z\in\C^+,\quad
1+zm_F(z)&=&\gamma+\gamma z m_{\underline{F}}(z),\label{eq:eq29}\\
\tau\left[1-\gamma^{-1}+\gamma^{-1}zm_F(z)\right]-z&=&-z\left[1+ m_{\underline{F}}(z)\tau\right].\label{eq:eq25}
\end{eqnarray}
Define: 
$$\mu(z)=\int_{-\infty}^{+\infty}\frac{g(\tau)}{1+ m_{\underline{F}}(z)\tau}\,dH(\tau).$$
Equation (\ref{eq:eq25}) yields:
\begin{eqnarray}
\lefteqn{\lim_{\varepsilon\to0^+}\lim_{\eta\to0^+}\frac{1}{\pi}\int_{-\varepsilon}^{+\varepsilon}
\int_{-\infty}^{+\infty}\im\left\{
\frac{\scriptstyle g(\tau)}{\scriptstyle\tau\left[1-\gamma^{-1}-\gamma^{-1}(\xi+i\eta)\,m_F(\xi+i\eta)\right]-\xi-i\eta}\right\}dH(\tau)d\xi}
\hspace{0.8cm}\nonumber\\
&=&\lim_{\varepsilon\to0^+}\lim_{\eta\to0^+}\frac{1}{\pi}\int_{-\varepsilon}^{+\varepsilon}
\im\left\{-\frac{1}{\xi+i\eta}
\int_{-\infty}^{+\infty}\frac{g(\tau)}{1+ m_{\underline{F}}(\xi+i\eta)\tau}\,dH(\tau)\right\}d\xi\nonumber\\
&=&\lim_{\varepsilon\to0^+}\lim_{\eta\to0^+}\frac{1}{\pi}\int_{-\varepsilon}^{+\varepsilon}
\im\left\{-\frac{\mu\left(\xi+i\eta\right)}{\xi+i\eta}\right\}\nonumber\\
&=&\int_{-\infty}^{+\infty}\frac{g(\tau)}{1+ \breve{m}_{\underline{F}}(0)\tau}\,dH(\tau),\nonumber
\end{eqnarray}
where the last equality follows from Lemma \ref{lem:zero}. $\square$
\paragraph{}
We are now ready to complete the proof of Theorem \ref{theo:kernel} for the case where $\gamma<1$. The inversion formula for the Stieltjes transform implies that:
\begin{eqnarray}
\lefteqn{\lim_{\varepsilon\to0^+}\left[\Phi(\varepsilon,\tau)-\Phi(-\varepsilon,\tau)\right]}\hspace{1cm}\nonumber\\
&=&\lim_{\varepsilon\to0^+}\lim_{\eta\to0^+}\frac{1}{\pi}\int_{-\varepsilon}^{+\varepsilon}
\im\left[\Theta^{\bone_{(-\infty,\tau)}}(\xi+i\eta)\right]d\xi\nonumber\\
&=&\lim_{\varepsilon\to0^+}\lim_{\eta\to0^+}\frac{1}{\pi}\int_{-\varepsilon}^{+\varepsilon}\im\left\{
\int_{-\infty}^{\tau}\frac{\scriptstyle dH(t)}{\scriptstyle t\left[1-\gamma^{-1}-\gamma^{-1}(\xi+i\eta)m_F(\xi+i\eta)\right]-\xi-i\eta}\right\}d\xi\nonumber\\
&=&\int_{-\infty}^{\tau}\frac{1}{1+\breve{m}_{\underline{F}}(0)\,t}\,dH(t),\label{eq:eq26}
\end{eqnarray}
where the last equality follows from Lemma \ref{lem:zero_g}. By Lemma \ref{lem:neighborhood}, we know that for $\lambda$ in a neighborhood of zero: $F(\lambda)=(1-\gamma)\bone_{[0,+\infty)}(\lambda)$. From Equation (\ref{eq:eq26}) we know that for $\lambda$ in a neighborhood of zero: 
$$\Phi(\lambda,\tau)=\int_{-\infty}^{\lambda}\int_{-\infty}^{\tau}
\frac{1}{1+\breve{m}_{\underline{F}}(0)\,t}\,dH(t)\,d\bone_{[0,+\infty)}(l).$$
Comparing the two expressions, we find that for $\lambda$ in a neighborhood of zero:
$$\Phi(\lambda,\tau)=\int_{-\infty}^{\lambda}\int_{-\infty}^{\tau}
\frac{1}{(1-\gamma)\left[1+\breve{m}_{\underline{F}}(0)\,t\right]}\,dH(t)\,dF(l).$$
Therefore, if we define $\varphi$ as in (\ref{eq:kernel}),
then we can see that  for $\lambda$ in a neighborhood of zero:
\be
\label{eq:phi}
\Phi(\lambda,\tau)=\int_{-\infty}^{\lambda}\int_{-\infty}^{\tau}\varphi(l,t)\,dH(t)\,dF(l).
\ee
From this point onwards, the fact that Equation (\ref{eq:phi}) holds for all $\lambda>0$ can be established exactly like we did in the case where $\gamma>1$. This completes the proof of Theorem \ref{theo:kernel}. $\square$

\section{Proofs of Theorems \ref{theo:Delta} and \ref{theo4b}\label{sec:remproofs}}
\subsection{Proof of Theorem \ref{theo:Delta} \label{sec:remproofs1}}

Lemma \ref{lem: convHi} shows that $\forall z\in\C^+,\,\Theta^{(1)}_N(z)\stackrel{\rm a.s.}{\to}\Theta^{(1)}(z)$, where:
\be
\label{eq:eq31}
\forall z\in\C^+,\, \Theta^{(1)}(z)
=\frac{\gamma}{1-\gamma^{-1}-\gamma^{-1}z\,m_F(z)}-\gamma.
\ee
Remember from Equation (\ref{eq:delta_link}) that:
$$\Delta_N(x)
=\lim_{\eta\to0^+}\frac{1}{\pi}\int_{-\infty}^x\im\left[\Theta_N^{(1)}(\lambda+i\eta)\right]d\lambda.$$
Therefore, by Lemma \ref{lem:omega}, $\lim_{N\to\infty}\Delta_N(x)$ exists and is equal to:
\be
\label{eq:eq32}
\Delta(x)=\lim_{\eta\to0^+}\frac{1}{\pi}\int_{-\infty}^x\im\left[\Theta^{(1)}(\lambda+i\eta)\right]d\lambda
\ee
for every $x\in\R$ where $\Delta$ is continuous. We first evaluate $\Delta(x)$ in the case where $\gamma>1$. Plugging Equation (\ref{eq:eq31}) into Equation (\ref{eq:eq32}) yields:
\begin{eqnarray}
\Delta(x)&=&\lim_{\eta\to0^+}\frac{1}{\pi}\int_{-\infty}^x\im\left[\frac{\gamma}{1-\gamma^{-1}-\gamma^{-1}(\lambda+i\eta)m_F(\lambda+i\eta)}-\gamma\right]d\lambda\nonumber\\
&=&\lim_{\eta\to0^+}\int_{-\infty}^x\frac{\pi^{-1}\im\left[(\lambda+i\eta)m_F(\lambda+i\eta)\right]}{\left|1-\gamma^{-1}-\gamma^{-1}(\lambda+i\eta)m_F(\lambda+i\eta)\right|^2}d\lambda\nonumber\\
&=&\int_{-\infty}^x\lim_{\eta\to0^+}\frac{\pi^{-1}\im\left[(\lambda+i\eta)m_F(\lambda+i\eta)\right]}{\left|1-\gamma^{-1}-\gamma^{-1}(\lambda+i\eta)m_F(\lambda+i\eta)\right|^2}d\lambda\label{eq:eq32b}\\
&=&\int_{-\infty}^x\frac{\pi^{-1}\im\left[\lambda\breve{m}_F(\lambda)\right]}{\left|1-\gamma^{-1}-\gamma^{-1}\lambda\breve{m}_F(\lambda)\right|^2}d\lambda\nonumber\\
&=&\int_{-\infty}^x\frac{\lambda F'(\lambda)}{\left|1-\gamma^{-1}-\gamma^{-1}\lambda\breve{m}_F(\lambda)\right|^2}d\lambda\nonumber\\
&=&\int_{-\infty}^x\frac{\lambda }{\left|1-\gamma^{-1}-\gamma^{-1}\lambda\,\breve{m}_F(\lambda)\right|^2}\,dF(\lambda),\nonumber
\end{eqnarray}
where Equation (\ref{eq:eq32b}) made use of Lemma \ref{lem:exchange}. This completes the proof of Theorem \ref{theo:Delta} in the case where $\gamma>1$.

In the case where $\gamma<1$, much of the arguments remain the same. The inversion formula for the Stieltjes transform implies that:
\begin{eqnarray}
\lefteqn{\lim_{\varepsilon\to0^+}\left[\Delta(\varepsilon)-\Delta(-\varepsilon)\right]}\hspace{1cm}\nonumber\\
&=&\lim_{\varepsilon\to0^+}\lim_{\eta\to0^+}\frac{1}{\pi}\int_{-\varepsilon}^{+\varepsilon}
\im\left[\Theta^{(1)}(\xi+i\eta)\right]d\xi\nonumber\\
&=&\lim_{\varepsilon\to0^+}\lim_{\eta\to0^+}\frac{1}{\pi}\int_{-\varepsilon}^{+\varepsilon}\im\left\{
\int_{-\infty}^{+\infty}\frac{\scriptstyle \tau\,dH(\tau)}{\scriptstyle \tau\left[1-\gamma^{-1}-\gamma^{-1}(\xi+i\eta)m_F(\xi+i\eta)\right]-\xi-i\eta}\right\}d\xi\nonumber\\
&=&\int_{-\infty}^{+\infty}\frac{\tau}{1+\breve{m}_{\underline{F}}(0)\,\tau}\,dH(\tau),\label{eq:eq33}
\end{eqnarray}
where the last equality follows from Lemma \ref{lem:zero_g}. Notice that for all $z\in\C^+$:
\begin{eqnarray}
\int_{-\infty}^{+\infty}\frac{\tau}{1+m_{\underline{F}}(z)\,\tau}\,dH(\tau)
&=&\frac{1}{m_{\underline{F}}(z)}\int_{-\infty}^{+\infty}\frac{1+m_{\underline{F}}(z)\,\tau-1}{1+m_{\underline{F}}(z)\,\tau}\,dH(\tau)\nonumber\\
&=&\frac{1}{m_{\underline{F}}(z)}
-\frac{1}{m_{\underline{F}}(z)}
\int_{-\infty}^{+\infty}\frac{1}{1+m_{\underline{F}}(z)\,\tau}\,dH(\tau).\label{eq:eq34}
\end{eqnarray}
Plugging Equation (\ref{eq:eq25}) into Equation (\ref{eq:eq34}) yields:
\begin{eqnarray}
\lefteqn{\int_{-\infty}^{+\infty}\frac{\tau}{1+m_{\underline{F}}(z)\,\tau}\,dH(\tau)}\hspace{2cm}\nonumber\\
&=&\frac{1}{m_{\underline{F}}(z)}
+\frac{z}{m_{\underline{F}}(z)}
\int_{-\infty}^{+\infty}\frac{1}{1-\gamma^{-1}+\gamma^{-1}zm_F(z)}\,dH(\tau)\nonumber\\
&=&\frac{1+z\,m_F(z)}{m_{\underline{F}}(z)},\label{eq:eq35}
\end{eqnarray}
where the last equality comes from the original Mar\v{c}enko-Pastur equation. Plugging Equation (\ref{eq:eq29}) into Equation (\ref{eq:eq35}) yields:
$$\int_{-\infty}^{+\infty}\frac{\tau}{1+m_{\underline{F}}(z)\,\tau}\,dH(\tau)
=\gamma\frac{1+z\,m_{\underline{F}}(z)}{m_{\underline{F}}(z)}.$$
Taking the limit as $z\in\C^+\to0$, we get:
$$\int_{-\infty}^{+\infty}\frac{\tau}{1+\breve{m}_{\underline{F}}(0)\,\tau}\,dH(\tau)
=\frac{\gamma}{\breve{m}_{\underline{F}}(0)}.$$
Plugging this result back into Equation (\ref{eq:eq33}) yields:
\be
\lim_{\varepsilon\to0^+}\left[\Delta(\varepsilon)-\Delta(-\varepsilon)\right]=\frac{\gamma}{\breve{m}_{\underline{F}}(0)}.\label{eq:eq36}
\ee
By Lemma \ref{lem:neighborhood}, we know that for $\lambda$ in a neighborhood of zero: $F(\lambda)=(1-\gamma)\bone_{[0,+\infty)}(\lambda)$. From Equation (\ref{eq:eq36}) we know that for $x$ in a neighborhood of zero: 
$$\Delta(x)=\int_{-\infty}^{x}
\frac{\gamma}{\breve{m}_{\underline{F}}(0)}\,d\bone_{[0,+\infty)}(\lambda).$$
Comparing the two expressions, we find that for $x$ in a neighborhood of zero:
$$\Delta(x)=\int_{-\infty}^{x}
\frac{\gamma}{(1-\gamma)\,\breve{m}_{\underline{F}}(0)}\,dF(\lambda).$$
Therefore, if we define $\delta$ as in (\ref{def: delta}),
then we can see that  for $x$ in a neighborhood of zero:
\be
\label{eq:eq37}
\Delta(x)=\int_{-\infty}^x\delta(\lambda)\,dF(\lambda).
\ee
From this point onwards, the fact that Equation (\ref{eq:eq37}) holds for all $x>0$ can be established exactly like we did in the case where $\gamma>1$. Thus the proof of Theorem \ref{theo:Delta} is complete. $\square$

\subsection{Proof of Theorem \ref{theo4b} \label{sec:remproofs2}}
As
\begin{eqnarray*}
&&\forall x\in\R,\quad\Psi_N(x)=\frac{1}{N}\sum_{i=1}^N\bone_{[\lambda_i,+\infty)}(x)
\sum_{j=1}^N\frac{\left|u_i^*v_j\right|^2}{\tau_j},\\
&&\forall z\in\C^+,\quad\Theta^{(-1)}_N(z)=\frac{1}{N}\sum_{i=1}^N\frac{1}{\lambda_i-z}
\sum_{j=1}^N\frac{\left|u_i^*v_j\right|^2}{\tau_j}, 
\end{eqnarray*}
and using the inversion formula for the Stieltjes transform, we obtain:
$$\forall x\in\R,\, \Psi_N(x)
=\lim_{\eta\to0^+}\int_{-\infty}^x\im\left[\Theta^{(-1)}_N(\lambda+i\eta)\right]d\lambda.$$
Since the function $g(\tau)=1/\tau$ is continuous on $[h_1,h_2]$, Theorem \ref{theo:theta} implies that $\forall z\in\C^+,$ $\Theta^{(-1)}_N(z)\stackrel{\rm a.s.}{\to}\Theta^{(-1)}(z)$, where:
\be
\label{eq:eq41}
\forall z\in\C^+,\quad\Theta^{(-1)}(z)
=\int_{-\infty}^{+\infty}\frac{\tau^{-1}}{\tau\left[1-\gamma^{-1}-\gamma^{-1}z\,m_F(z)\right]-z}\,dH(\tau).
\ee 
Therefore, by Lemma \ref{lem:omega}, $\lim_{N\to\infty}\Psi_N(x)$ exists and is equal to:
\be
\label{eq:eq42}
\Psi(x)=\lim_{\eta\to0^+}\frac{1}{\pi}\int_{-\infty}^x\im\left[\Theta^{(-1)}(\lambda+i\eta)\right]d\lambda,
\ee
for every $x\in\R$ where $\Psi$ is continuous. We first evaluate $\Psi(x)$ in the case where $\gamma>1$, so that $F$ is continuously differentiable on all of $\R$. 

In the notation of Lemma \ref{lem:continuous}, we set $\nu$ equal to zero so that $\forall\tau\in\R, \quad\rho(\tau)=g(\tau)\times \tau =1$. Then Equation (\ref{eq:eq9}) implies that:
$$\forall z\in\C^+,\quad\Theta^{(-1)}(z)
=\frac{m_F(z)-\left[1+\gamma^{-1}\Theta^{(1)}(z)\right]\int_{-\infty}^{+\infty}\tau^{-1}dH(\tau)}{
z\left[1+\gamma^{-1}\Theta^{(1)}(z)\right]}.$$
Using Equation (\ref{eq:eq4}), we obtain:
\be
\label{eq:eq43}
\forall z\in\C^+,\quad\Theta^{(-1)}(z)
=\frac{m_F(z)}{z}\left[1-\gamma^{-1}-\gamma^{-1}z\,m_F(z)\right]
-\frac{1}{z}\int_{-\infty}^{+\infty}\tau^{-1}dH(\tau).
\ee
Thus for all $\lambda\in\R$:
\begin{eqnarray*}
\lim_{\eta\to0^+}
\im\left[\Theta^{-1}(\lambda+i\eta)\right]
&=&\frac{1}{\lambda}\,\im\left\{\breve{m}_F(\lambda)\left[1-\gamma^{-1}-\gamma^{-1}\lambda\,\breve{m}_F(\lambda)\right]\right\}\\
&=&\frac{1}{\lambda}\left\{1-\gamma^{-1}-2\gamma^{-1}\lambda\,\re\left[\breve{m}_F(\lambda)\right]\right\}
\times\im\left[\breve{m}_F(\lambda)\right]\\
&=&\frac{1}{\lambda}\left\{1-\gamma^{-1}-2\gamma^{-1}\lambda\,\re\left[\breve{m}_F(\lambda)\right]\right\}
\times\pi F'(\lambda).
\end{eqnarray*}
Plugging this result back into Equation (\ref{eq:eq42}) yields:
\begin{eqnarray*}
\Psi(x)&=&\frac{1}{\pi}\int_{-\infty}^x\lim_{\eta\to0^+}\im\left[\Theta^{(-1)}(\lambda+i\eta)\right]d\lambda\\
&=&\int_{-\infty}^x\frac{1-\gamma^{-1}-2\gamma^{-1}\lambda\,\re\left[\breve{m}_F(\lambda)\right]}{\lambda}
\,dF(\lambda),
\end{eqnarray*}
where we made use of Lemma \ref{lem:exchange}. This completes the proof of Theorem \ref{theo4b} in the case where $\gamma>1$.

We now turn to the case where $\gamma<1$. Equation (\ref{eq:eq42}) implies that:
\be
\label{eq:eq44}
\lim_{\varepsilon\to0^+}\left[\Psi(\varepsilon)-\Psi(-\varepsilon)\right]
=\lim_{\varepsilon\to0^+}\lim_{\eta\to0^+}\frac{1}{\pi}\int_{-\varepsilon}^{+\varepsilon}
\im\left[\Theta^{(-1)}(\xi+i\eta)\right]d\xi.
\ee
Plugging Equation (\ref{eq:eq29}) into Equation (\ref{eq:eq43}) yields for all $z\in\C^+$:
\begin{eqnarray*}
\Theta^{(-1)}(z)
&=&-m_F(z)\,m_{\underline{F}}(z)-\frac{1}{z}\int_{-\infty}^{+\infty}\frac{1}{\tau-0}dH(\tau)\\
&=&\frac{1}{z}\left[1-\gamma-\gamma zm_{\underline{F}}(z)\right]
m_{\underline{F}}(z)-\frac{1}{z}\breve{m}_H(0).
\end{eqnarray*}
Plugging this result into Equation (\ref{eq:eq44}), we get:
$$\lim_{\varepsilon\to0^+}\left[\Psi(\varepsilon)-\Psi(-\varepsilon)\right]
=\lim_{\varepsilon\to0^+}\lim_{\eta\to0^+}\frac{1}{\pi}\int_{-\varepsilon}^{+\varepsilon}
\im\left\{-\frac{\mu(\xi+i\eta)}{\xi+i\eta}
\right\}d\xi,$$
where $\mu(z)=-[1-\gamma-\gamma zm_{\underline{F}}(z)]m_{\underline{F}}(z)+\breve{m}_H(0)$. Therefore, by Lemma \ref{lem:zero}, we have:
\be
\label{eq:eq45}
\lim_{\varepsilon\to0^+}\left[\Psi(\varepsilon)-\Psi(-\varepsilon)\right]=\mu(0)
=-(1-\gamma)\breve{m}_{\underline{F}}(0)+\breve{m}_H(0).
\ee
By Lemma \ref{lem:neighborhood}, we know that for $\lambda$ in a neighborhood of zero: $F(\lambda)=(1-\gamma)\bone_{[0,+\infty)}(\lambda)$. From Equation (\ref{eq:eq45}) we know that for $x$ in a neighborhood of zero: 
$$\Psi(x)=\int_{-\infty}^{x}
\left[-(1-\gamma)\breve{m}_{\underline{F}}(0)+\breve{m}_H(0)\right]
\,d\bone_{[0,+\infty)}(\lambda).$$
Comparing the two expressions, we find that for $x$ in a neighborhood of zero:
$$\Psi(x)=\int_{-\infty}^{x}
\left[-\breve{m}_{\underline{F}}(0)+\frac{1}{1-\gamma}\,\breve{m}_H(0)\right]dF(\lambda).$$
Therefore, if we define $\psi$ as in (\ref{def: psi}),
then we can see that for $x$ in a neighborhood of zero:
\be
\label{eq:eq46}
\Psi(x)=\int_{-\infty}^x\psi(\lambda)\,dF(\lambda).
\ee
From this point onwards, the fact that Equation (\ref{eq:eq46}) holds for all $x>0$ can be established exactly like we did in the case where $\gamma>1$. Thus the proof of Theorem \ref{theo4b} is complete. $\square$

\paragraph{Acknowledgements}
O. Ledoit wishes to thank the organizers and participants of the Stanford Institute for Theoretical Economics (SITE) summer 2008 workshop on ``Complex Data in Economics and Finance'' for their comments on an earlier version of this paper. S. P\'ech\'e thanks Prof. J. Silverstein for his helpful comments on a preliminary version of this paper.

\bibliographystyle{spmpsci}
\bibliography{biblio1}


\end{document}